\documentclass[12pt]{article}
\usepackage{amstex}
\usepackage{amssymb}
\usepackage{amscd}
\usepackage{a4}
\usepackage{pb-diagram}

\def\bk{\mathbf k\ }
\def\gf{{\hbox to 12mm{\rightarrowfill}}}
\def\df{{\hbox to 12mm{\downarrowfill}}}

\def\oisi{\buildrel\textstyle{\rm i \atop ^\vee} \over\dots}
\def\omis{\buildrel {\rm i \atop ^\vee} \over .}

\def\Talpha#1{\vbox{\ialign{##\crcr
$\alpha$\crcr\noalign{\kern2pt\nointerlineskip}
$\hfil\displaystyle{#1}\hfil$\crcr}}} 
 
\def\Onabla#1{\vbox{\ialign{##\crcr$\,\scriptstyle{0}$\crcr
\noalign{\kern2pt\nointerlineskip}
$\hfil\displaystyle{#1}\hfil$\crcr}}}
\def\im{\mbox{Im}}
\def\mod{\mbox{mod}}

 \def\cale{{\cal E}}
 \def\cala{{\cal A}}
\def\fraca{{\mathfrak{A}}}

\def\calm{{\cal M}}
\def\caln{{\cal N}}
\def\cals{{\cal S}}
\def\calh{{\cal H}}

\def\fg{\mathfrak f}
\def\Fg{\mathfrak F}
\def\ft{\mathfrak T}
\def\sg{\mathfrak s}
\def\Sg{\mathfrak S}
\def\Bg{\mathfrak B}
 \def\oplusinf{\mathop{\oplus}}
 \def\cals{{\cal S}}

\def\cals{{\cal S}}
 
 \def\calf{{\cal F}}
 
\def\bbbone{\mbox{\rm 1\hspace {-.6em} l}}

\def\gder{{\mbox{Der}}}

\def\fil{\displaystyle{\buildrel {[i]^\ell}\over \longrightarrow}}
\def\fdm{\displaystyle{\buildrel {[d]^m}\over \longrightarrow}}

\newtheorem{corollary}{COROLLARY}
 \newtheorem{theorem}{THEOREM}
\newtheorem{lemma}{LEMMA} 
\newtheorem{proposition}{PROPOSITION} 
\begin{document}

\baselineskip=0.9cm
\vspace{2cm}

\begin{center} {\Huge\bf $d^N=0$} \end{center} \vspace{0.75cm}

\begin{center} Michel DUBOIS-VIOLETTE \\ \vspace{0.3cm} {\small Laboratoire de
Physique Th\'eorique et Hautes Energies\footnote{Laboratoire associ\'e au
Centre
National de la Recherche Scientifique - URA D0063}\\ Universit\'e Paris XI,
B\^atiment 211\\ 91 405 Orsay Cedex, France\\ flad$@$qcd.th.u-psud.fr}

\end{center} \vspace{1cm}

\begin{center} \today \end{center}

\vspace {1cm}

\begin{abstract}
We study the generalized homology associated with a nilpotent endomorphism $d$
satisfying $d^N=0$. For simplicial modules we construct such nilpotent
endomorphisms and we prove a general result relating the corresponding
generalized homologies to the ordinary homology. We also discuss the
generalization of the notion of graded differential algebra in this context.
\end{abstract}

\vspace {5cm}

\noindent L.P.T.H.E.-ORSAY 97/53\\
\noindent q-alg/9710021

\newpage

\section{Introduction}

The aim of this work is to describe and study the generalized homologies
associated with nilpotent endomorphisms $d$ satisfying $d^N=0$, $N$ being an
integer with $N\geq 2$. Throughout this paper $\bk$ is a commutative ring with
a unit and $N$ is an integer greater than or equal to 2.\\

In Section 2 the notion of $N$-differential module is introduced and the
corresponding generalized homology is defined and studied. A $N$-differential
module over $\bk$ is a $\bk$-module equipped with an endomorphism $d$, its
$N$-differential, satisfying $d^N=0$; its homology being the $\bk$-modules
$\ker (d^n)/\im (d^{N-n})$, $n=1,\dots,N-1$. Lemma 1, which has no nontrivial
counterpart in the classical situation $N=2$, turns out to be a basic result
for the development of the theory.\\

In Section 3 the graded case is studied. A $N$-complex of modules is a
$N$-differential module which is graded with a $N$-differential homogeneous of
degree 1 or $-1$. To a $N$-complex is associated a family of ordinary complexes
and we study the relations between their homologies. Given a pre-(co)simplicial
module and an element $q$ of the ring $\bk$ such that $[N]_q=0$, (see in
Section 2 for the definitions of ``$q$-numbers"), we construct a family of
$N$-differentials and we show that the corresponding generalized (co)homologies
are related by homomorphisms to the ordinary (co)homology of the
pre-(co)simplicial module. In the case of a (co)simplicial module and if
furthermore $q$ is such that $[n]_q$ is invertible for $n=1,\dots,N-1$ we show
that these homomorphisms are isomorphisms and that this allows to identify
these generalized (co)homologies in terms of the ordinary one. Although that
kind of results may be expected in view of the Dold-Kan correspondence, the
identifications as well as the proof of these identifications are far from
being trivial. Applied to the generalized Hochschild (co)homology at roots of
unity, these results imply the results announced in \cite{MD-V} as well as the
results of \cite{KW}. Applied to simplicial complexes in the case where $N$ is
a prime number and where $\bk=\mathbb Z/N\mathbb Z$, these results give an
improvement of the results of \cite{May}.\\

In Section 4 the corresponding generalization of graded differential algebras
is introduced and discussed. Several examples are described and we construct
universal models associated with associative algebras generalizing thereby the
construction of the universal graded differential envelope of an associative
algebra \cite{Kar}. These examples have been introduced in \cite{D-VK} or in
\cite{MD-V} in the case where $\bk=\mathbb C$, however here we give a different
presentation in the framework of cosimplicial algebras. Furthermore there are
some subtilities in the replacement of the field $\mathbb C$ by the ring $\bk$.

\newpage

\section{$N$-differential modules}

Throughout this paper $\bk$ is a commutative ring with a unit and $N$ denotes
an integer with $N\geq 2$. By a module without other specification we always
mean a $\bk$-module; the same convention is adopted for endomorphisms,
homomorphisms, etc.. A $N$-{\sl differential module} is a module $E$ equipped
with an endomorphism $d$, its $N$-{\sl differential}, satisfying $d^N=0$. Given
two $N$-differential modules $(E,d)$ and $(E',d')$, an {\sl homomorphism of
$N$-differential modules} of $E$ into $E'$ is an homomorphism
$\varphi:E\rightarrow E'$ satisfying $\varphi\circ d=d'\circ \varphi$. The
$N$-differential modules and their homomorphisms form an abelian category so
one can speak there of kernels, exact sequences, etc.\\

Let $E$ be a $N$-differential module with $N$-differential $d$. One defines for
each integer $n$ with $1\leq n\leq N-1$ the submodules $Z_{(n)}(E)$ and
$B_{(n)}(E)$ by $Z_{(n)}(E)=\ker (d^n)$ and $B_{(n)}(E)= \mathrm{Im\ }
(d^{N-n})$. Since $d^N=0$, one has $B_{(n)}(E)\subset Z_{(n)}(E)$ so the module
$H_{(n)}(E)=Z_{(n)}(E)/B_{(n)}(E)$ is well defined. The modules $H_{(n)}(E),\
1\leq n\leq N-1$, will be refered to as {\sl the homology of} $E$. Assume that
$N\geq 3$ and let $n$ be an integer with $1\leq n\leq N-2$, one has
$Z_{(n)}(E)\subset Z_{(n+1)}(E)$ and $B_{(n)}(E)\subset B_{(n+1)}(E)$ which
induces an homomorphism $[i]:H_{(n)}(E)\rightarrow H_{(n+1)}(E)$; on the other
hand one also has $dZ_{(n+1)}(E)\subset Z_{(n)}(E)$ and $dB_{(n+1)}(E)\subset
B_{(n)}(E)$ which induces an homomorphism $[d]:H_{(n+1)}(E)\rightarrow
H_{(n)}(E)$. The following basic lemma \cite{D-VK} shows that the $H_{(n)}(E)$
are not independent.
\newpage

\begin{lemma}

Let $\ell$ and $m$ be integers with $\ell\geq 1,\ m\geq 1$ and $\ell+m\leq
N-1$. Then the following hexagon $(\calh^{\ell,m})$ of homomorphisms
\[ \begin{diagram} \node{}
\node{H_{(\ell+m)}(E)} \arrow{e,t}{[d]^m} \node{H_{(\ell)}(E)}
\arrow{se,t}{[i]^{N-(\ell+m)}} \node{} \\ \node{H_{(m)}(E)}
\arrow{ne,t}{[i]^\ell}
\node{} \node{} \node{H_{(N-m)}(E)} \arrow{sw,b}{[d]^\ell} \\ \node{}
\node[1]{H_{(N-\ell)}(E)} \arrow{nw,b}{[d]^{N-(\ell+m)}}
\node{H_{(N-(\ell+m))}(E)}
\arrow{w,b}{[i]^m} \node{} \end{diagram} \]
is exact.
\end{lemma}
\noindent
{\bf Proof}. \cite{D-VK}. By interchanging $(m,\ell),(\ell,N-(\ell+m))$ and
$(N-(\ell+m),m)$, it is sufficient proving the exactness of the sequences
$H_{(m)}(E)\fil H_{(\ell+m)}(E)\fdm H_{(\ell)}(E)$ and $H_{(\ell+m)}(E) \fdm
H_{(\ell)}(E)
\mathop{\hbox to 12mm{\rightarrowfill}}^{[i]^{N-(\ell+m)}} H_{(N-m)}(E)$. It
follows from the definitions that $[d]^m\circ[i]^\ell$ is the zero map of
$H_{(m)}(E)$ into $H_{(\ell)}(E)$ and that $[i]^{N-(\ell+m)}\circ [d]^m$ is the
zero map of $H_{(\ell+m)}(E)$ into $H_{(N-m)}(E)$. Let $z\in Z_{(\ell+m)}(E)$
be such that $d^mz=d^{N-\ell}c$ for some $c\in E$; then
$d^m(z-d^{N-(\ell+m)}c)=0$ which means that the class of $z$ in
$H_{(\ell+m)}(E)$ belongs to $[i]^\ell H_{(m)}(E)$ which achieves the proof of
the exactness of the first sequence. Let $z\in Z_{(\ell)}(E)$ be such that
$z=[d]^mc$ for some $c\in E$; then $d^{\ell+m}c=0$ which means that the class
of $z$ in $H_{(\ell)}(E)$ belongs to $[d]^mH_{(\ell+m)}(E)$ which achieves the
proof of the exactness of the second sequence. $\Box$\\

Let $\varphi:E\rightarrow E'$ be an homomorphism of $N$-differential modules.
Then one has $\varphi(Z_{(n)}(E))\subset Z_{(n)}(E')$ and
$\varphi(B_{(n)}(E))\subset B_{(n)}(E')$ so $\varphi$ induces an homomorphism
$\varphi_\ast:H_{(n)}(E)\rightarrow H_{(n)}(E')$ for each integer $n$ with
$1\leq n\leq N-1$. The correspondences $E\mapsto H_{(n)}(E)$ are of course
functorial with $H_{(n)}(\varphi)=\varphi_\ast$. Furthermore one has
$\varphi_\ast\circ [i]=[i]\circ \varphi_\ast$ and $\varphi_\ast\circ
[d]=[d]\circ \varphi_\ast$.

\begin{proposition}
Let $\varphi:E\rightarrow E'$ be an homomorphism of $N$-differential modules
and assume that $\varphi_\ast:H_{(1)}(E)\rightarrow H_{(1)}(E')$ and
$\varphi_\ast:H_{(N-1)}(E)\rightarrow H_{(N-1)}(E')$ are isomorphisms. Then
$\varphi_\ast:H_{(n)}(E)\rightarrow H_{(n)}(E')$ is an isomorphism for any
integer $n$ with $1\leq n\leq N-1$.
\end{proposition}

\noindent{\bf Proof.} By induction: Let $n$ be an integer with $1\leq n\leq
N-2$ and assume that $\varphi_\ast:H_{(m)}(E)\rightarrow H_{(m)}(E')$ and
$\varphi_\ast:H_{(N-m)}(E)\rightarrow H_{(N-m)}(E')$ are isomorphisms for each
integer $m$ such that $1\leq m\leq n$; considering the action of $\varphi_\ast$
on the hexagon $\calh^{n,1}$, one obtains the two commutative diagrams with
exact columns
\[
\begin{diagram}
\node{H_{(N-n)}(E)}\arrow{e,t}{\varphi_\ast}\arrow{s,r}{[d]^{N-n-1}}
\node{H_{(N-n)}(E')}\arrow{s,r}{[d]^{N-n-1}}\arrow{e,!}
\node{H_{(n)}(E)}\arrow{e,t}{\varphi_\ast}\arrow{s,r}{[i]^{N-n-1}}
\node{H_{(n)}(E')}\arrow{s,r}{[i]^{N-n-1}}\\
\node{H_{(1)}(E)}\arrow{e,t}{\varphi_\ast}\arrow{s,r}{[i]^n}
\node{H_{(1)}(E')}\arrow{s,r}{[i]^n}\arrow{e,!}
\node{H_{(N-1)}(E)}\arrow{e,t}{\varphi_\ast}\arrow{s,r}{[d]^n}
\node{H_{(N-1)}(E')}\arrow{s,r}{[d]^n}\\
\node{H_{(n+1)}(E)}\arrow{e,t}{\varphi_\ast}\arrow{s,r}{[d]}
\node{H_{(n+1)}(E')}\arrow{s,r}{[d]}\arrow{e,!}
\node{H_{(N-n-1)}(E)}\arrow{e,t}{\varphi_\ast}\arrow{s,r}{[i]}
\node{H_{(N-n-1)}(E')}\arrow{s,r}{[i]}\\
\node{H_{(n)}(E)}\arrow{e,t}{\varphi_\ast}\arrow{s,r}{[i]^{N-n-1}}
\node{H_{(n)}(E')}\arrow{s,r}{[i]^{N-n-1}}\arrow{e,!}
\node{H_{(N-n)}(E)}\arrow{e,t}{\varphi_\ast}\arrow{s,r}{[d]^{N-n-1}}
\node{H_{(N-n)}(E')}\arrow{s,r}{[d]^{N-n-1}}\\
\node{H_{(N-1)}(E)}\arrow{e,t}{\varphi_\ast}
\node{H_{(N-1)}(E')}\arrow{e,!}
\node{H_{(1)}(E)}\arrow{e,t}{\varphi_\ast}
\node{H_{(1)}(E')}
\end{diagram}
\]
in which the first two rows and the last two rows are isomorphisms, which
implies by the ``five lemma" that $\varphi_\ast:H_{(n+1)}(E)\rightarrow
H_{(n+1)}(E')$ and $\varphi_\ast: H_{(N-(n+1))}(E)\rightarrow
H_{(N-(n+1))}(E')$ are also isomorphisms. $\Box$\\

\noindent {\bf Remark 1.} It is worth noticing that $H_{(1)}(E)\oplus
H_{(N-1)}(E)$ is the ordinary homology of the ordinary ($\mathbb Z_2$-graded)
differential module $\cale=E\oplus E$ with differential $\delta$ defined by
$\delta(x\oplus y)=d^{N-1}(y)\oplus d(x)$.\\

Let us now generalize in this context the construction of the connecting
homomorphism associated to a short exact sequence of differential modules. So
let $0\rightarrow E \stackrel{\alpha}{\rightarrow} F
\stackrel{\beta}{\rightarrow}G\rightarrow 0$ be a short exact sequence of
$N$-differential modules and let $g\in Z_{(n)}(G),\ (1\leq n\leq N-1)$. Then,
by exactness at $G$, there is an element $f$ of $F$ with $\beta(f)=g$. One has
$\beta(d^nf)=d^ng=0$ and therefore, by exactness at $F$, there is an element
$e$ of $E$ with $\alpha(e)=d^n f$. One has $\alpha(d^{N-n}e)=d^Nf=0$ which
implies, by exactness at $E$, $d^{N-n}e=0$ that is $e\in Z_{(N-n)}(E)$. Given
$g$ as above, $f$ is unique up to $\alpha(E)$, in view of exactness at $F$, and
therefore $e$ is unique up to $d^nE=B_{(N-n)}(E)$ in view of exactness of $E$.
Furthermore if $g\in d^{N-n}G=B_{(n)}(G)$ then $e\in d^nE=B_{(N-n)}(E)$.
Therefore one has a well defined homomorphism $\partial:H_{(n)}(G)\rightarrow
H_{(N-n)}(E)$ which generalizes the usual connecting homomorphism. These
homomorphisms $\partial:H_{(n)}(G)\rightarrow H_{(N-n)}(E)$, for $1\leq n\leq
N-1$, will be refered to as the {\sl connecting homomorphisms} associated with
the short exact sequence $0\rightarrow E \stackrel{\alpha}{\rightarrow} F
\stackrel{\beta}{\rightarrow}G\rightarrow 0$ of $N$-differential modules. This
construction is functorial in the following sense: If one has a commutative
diagram of $N$-differential modules
\[
\begin{diagram}
\node{0} \arrow{e} \node{E} \arrow{e,t}{\alpha} \arrow{s,l}{\lambda}
\node{F}\arrow{e,t}{\beta} \arrow{s,l}{\mu} \node{G} \arrow{e} \arrow{s,l}{\nu}
\node{0}\\
\node{0} \arrow{e} \node{E'} \arrow{e,b}{\alpha'} \node{F'}
\arrow{e,b}{\beta'}\node{G'}\arrow{e}\node{0}
\end{diagram}
\]
where the rows are exact, then the connecting homomorphisms $\partial$ and
$\partial'$ associated with the first and the second row are related by
$\lambda_\ast\circ\partial=\partial'\circ \nu_\ast: H_{(n)}(G)\rightarrow
H_{(N-n)}(E')$.\\

\noindent {\bf Remark 2.} For later purpose let us notice the following:\\
For each $n$ ($1\leq n\leq N-1$), the application of $d^n$ enter in the
construction of $\partial:H_{(n)}(G)\rightarrow H_{(N-n)}(E)$ in such a way
that, in the graded situation, the degree of $\partial:H_{(n)}(G)\rightarrow
H_{(N-n)}(E)$ is the degree of $d^n$.\\

The following lemma is the generalization to $N$-modules of the exact triangle
in homology associated with a short exact sequence of differential modules (the
usual case $N=2$).

\begin{lemma}
Let $0\rightarrow E
\stackrel{\alpha}{\rightarrow}F\stackrel{\beta}{\rightarrow} G \rightarrow 0$
be a short exact sequence of\linebreak[4] $N$-differential modules with
associated connecting homomorphisms denoted by $\partial$. Then for each
integer $n$ with $1\leq n\leq N-1$, the following hexagon $(\calh_n)$ of
homomorphisms

\[ \begin{diagram} \node{}
\node{H_{(n)}(F)} \arrow{e,t}{\beta_\ast} \node{H_{(n)}(G)}
\arrow{se,t}{\partial} \node{} \\ \node{H_{(n)}(E)} \arrow{ne,t}{\alpha_\ast}
\node{} \node{} \node{H_{(N-n)}(E)} \arrow{sw,b}{\alpha_\ast} \\ \node{}
\node[1]{H_{(N-n)}(G)} \arrow{nw,b}{\partial} \node{H_{(N-n)}(F)}
\arrow{w,b}{\beta_\ast} \node{} \end{diagram} \]

is exact.
\end{lemma}

The proof of this lemma follows from the definition and is the same as the
proof of exactness of the triangle in the usual case $N=2$. Kassel and Wambst
gave in [KW] a direct proof of this lemma together with the construction of
$\partial$ by using the snake lemma.\\

Let us give now some criteria ensuring the vanishing of the $H_{(n)}(E)$. The
first criterion is extracted from \cite{Kapr}.

\begin{lemma}
Let $E$ be a $N$-differential module such that $H_{(k)}(E)=0$ for some integer
$k$ with $1\leq k\leq N-1$. Then one has $H_{(n)}(E)=0$ for any integer $n$
with $1\leq n\leq N-1$.
\end{lemma}

\noindent {\bf Proof}. Although Kapranov gave a nice direct proof of this in
\cite{Kapr}, we give here a proof by using  Lemma 1. One extracts from the
hexagon $(\calh^{k-1,1})$, when $k\geq 2$, the exact sequence
\[
H_{(N-1)}(E)\mathop{\hbox to 12mm{\rightarrowfill}}^{[d]^{N-k}}
H_{(k-1)}(E)\stackrel{[i]}{\rightarrow} H_{(k)}(E)=0
\]
and the mapping $[d]^{N-k}:H_{(N-1)}(E)\rightarrow H_{(k-1)}(E)$ is the zero
mapping because either $k=N-1$ or if $k\leq N-2$ it factorizes through
\[
H_{(N-1)}(E)\mathop{\hbox to
12mm{\rightarrowfill}}^{[d]^{N-k-1}}H_{(k)}(E)=0\stackrel{[d]}{\rightarrow}H_{(k-1)}(E).
\]
Thus $H_{(k)}(E)=0$ implies $H_{(n)}(E)=0$ for any $n$ with $1\leq n\leq k$. On
the other hand, one extracts from the hexagon $(\calh^{k,1})$, when $k\leq
N-2$, the exact sequence $0=H_{(k)}(E) \stackrel{[i]}{\rightarrow} H_{(k+1)}(E)
\stackrel{[d]^k}{\rightarrow} H_{(1)}(E)$ and $H_{(1)}(E)=0$ in view of above,
therefore $H_{(k+1)}(E)=0$. This achieves the proof of the lemma.$\Box$\\

For next criterion which is in \cite{Kapr} and in \cite{KW}, we introduce,
following Kassel and Wambst \cite{KW}, the appropriate notion of homotopy. Let
$E$ and $F$ be two $N$-differential modules; two homomorphisms of
$N$-differential modules, $\lambda$ and $\mu$, of $E$ into $F$ will be said to
be {\sl homotopic} if there exist module-homomorphisms $h_k:E\rightarrow F$,
$(k=0,1,\dots, N-1)$, such that $\lambda-\mu=\displaystyle{\sum^{N-1}_{k=0}}
d^{N-1-k}h_k d^k$. This is an equivalence relation.
\begin{lemma}
If $\lambda$ and $\mu$ are homotopic homomorphisms of $N$-differential modules
of $E$ into $F$, they induce the same homomorphisms in homology, i.e.
\[
\lambda_\ast=\mu_\ast: H_{(n)}(E)\rightarrow H_{(n)}(F)
\]
for each integer $n$ with $1\leq n\leq N-1$.
\end{lemma}

\noindent {\bf Proof}. \cite{KW}. By applying
$\lambda-\mu=\displaystyle{\sum^{N-1}_{k=0}} d^{N-1-k}h_kd^k$ to an element $z$
of $Z_{(n)}(E)$ one obtains
$\lambda(z)-\mu(z)=d^{N-n}\left(\displaystyle{\sum^{n-1}_{k=0}d^{n-1-k}h_kd^kz)}\right),\ (1\leq n\leq N-1)$. Thus the difference between the elements $\lambda(z)$ and $\mu(z)$ of $Z_{(n)}(F)$ belongs to $B_{(n)}(F)$ which implies the statement. $\Box$\\
Next criterion follows immediately.
\begin{corollary}
Let $E$ be an $N$-differential module such that there are endomorphisms of
modules $h_k:E\rightarrow E$ for $k=0,1,\dots,N-1$ satisfying \linebreak[4]
$\displaystyle{\sum^{N-1}_{k=0}}d^{N-1-k}h_kd^k=Id_E$; then one has
$H_{(n)}(E)=0$ for each integer $n$ with\linebreak[4] $1\leq n\leq N-1$.
\end{corollary}

In order to formulate the last criterion which will be the most useful in the
following, we recall the definition of $q$-numbers. Let $q$ be an element of
the ring $\bk$, one associates with $q$ a mapping $[.]_q:\mathbb N \rightarrow
\bk$, $n\mapsto [n]_q$,  which is defined by setting
$[0]_q=0$ and $[n]_q=1+\dots+q^{n-1}=\sum^{n-1}_{k=0} q^k$
for $n\geq 1,\ (q^0=1)$. For $n\in \mathbb N$ with $n\geq 1$, one defines the
$q$-factorial $[n]_q!\in \bk$ by $[n]_q\dots 1=\prod^n_{k=1}[k]_q$. For
integers $n$ and $m$ with $n\geq 1$ and $0\leq m\leq n$, one defined
inductively the $q$-binomial coefficients $\left[\begin{array}{c} n\\m
\end{array}\right]_q\in \bk$ by setting $\left[\begin{array}{c} n\\0
\end{array}\right]_q
= \left[\begin{array}{c} n\\n \end{array}\right]_q=1$ and
$\left[\begin{array}{c} n\\m \end{array}\right]_q +
q^{m+1}\left[\begin{array}{c} n\\m+1
\end{array}\right]_q=\left[\begin{array}{c} n+1\\m+1 \end{array}\right]_q$ for
$0\leq m\leq n-1$. As in \cite{KW} let us introduce the following assumptions
$(A_0)$ and $(A_1)$ on the ring $\bk$ and the element $q$ of $\bk$ :\\
$(A_0)$\hspace{1cm} $[N]_q=0$\\
$(A_1)$\hspace{1cm} $[N]_q=0$ and $[n]_q$ is invertible for $1\leq n\leq N-1$,
$(n\in \mathbb N)$.\\
Notice that $[N]_q=0$ implies that $q^N=1$ and therefore that $q$ is
invertible. Furthermore if $q$ is invertible one has
$[n]_{q^{-1}}=q^{-n+1}[n]_q$, $\forall n\in \mathbb N$. Therefore Assumption
$(A_0)$, (resp. $(A_1)$), for $\bk$ and $q\in \bk$ is equivalent to Assumption
$(A_0)$, (resp. $(A_1)$), for $\bk$ and $q^{-1}\in \bk$.\\
We are now ready to state the last criterion \cite{MD-V}.

\begin{lemma}
Suppose that $\bk$ and $q\in\bk$ satisfy $(A_1)$ and let $E$ be a
$N$-differential module. Assume that there is a module-endomorphism $h$ of $E$
such that\linebreak[4] $hd-qdh=Id_E$. Then one has $H_{(n)}(E)=0$ for each
integer $n$ with $1\leq n\leq N-1$.
\end{lemma}

\noindent {\bf Proof}. We shall show that one has the identity
$\displaystyle{\sum^{N-1}_{k=0}}d^{N-1-k}h^{N-1} d^k=[N-1]_q! Id_E$ which
implies the result in view of Corollary 1. In order to prove the latter
identity, consider the unital associative $\bk$-algebra ${\cal D}_q$ generated
by two elements $H$ and $D$ with relation $HD-qDH=\bbbone$. In ${\cal D}_q$,
one clearly has $\displaystyle{\sum^{N-1}_{k=0}} D^{N-1-k} H^{N-1}
D^k=\displaystyle{\sum^{N-1}_{r=0}} a_r D^r H^r$ for some $a_r\in \bk$. We
shall show that, under assumption $(A_1)$, one has $a_0=[N-1]_q!$ and $a_r=0$
for $1\leq r\leq N-1$ which implies the result. There is a proof of this in
\cite{KW}. Here however we give a simpler proof which has been suggested by
Rausch de Traubenberg and which is essentially extracted from the algebraic
technics of \cite{RdT}. Consider the matrices $H_N$ and $D_N$ of $M_N(\bk)$
defined by

\[
H_N=\left(
\begin{array}{ccccccc}
0 & . & . & . & . & . & 0\\
{[N-1]}_q & . & & & & &.\\
0 & . & . &&&& .\\
. & . & . & . &&& .\\
. &  & . & . & . & &.\\
. & & &.&.&.&.\\
0 & . &. &. & 0 & [1]_q & 0\\
\end{array}
\right)
\]
 \mbox{and}
\[
D_N=\left(
\begin{array}{ccccccc}
0 & 1 & 0 & . & . & . & 0\\
. & . & . & . &  & &.\\
. &   & . & . & . &  & .\\
. &   &   & . & . & . & .\\
. &  &   &   & .  & . & 0\\
. & & & & &.& 1\\
0 & . &. &. & . & . & 0\\
\end{array}
\right),
\]
$([1]_q=1).$\\
One has in $M_N(\bk)$ :

\[
H_ND_N-qD_NH_N=\left(
\begin{array}{ccccccc}
-q[N-1]_q &   &   &   &   &   &  \\
 & 1 &   &   &   &   &  \\
  &   & .  &   &   &   & \\
  &   &   & . &   &   &  \\
  &  &   &   & .  &   &  \\
  & & & & &.&  \\
  &   &  &  &   &   & 1\\
\end{array}
\right).
\]
It follows from $[N]_q=0$ that $-q[N-1]_q=1$ and thus $H_ND_N-qD_NH_N=\bbbone$,
which means that $H\mapsto H_N$, $D\mapsto D_N$ define an homomorphism of
${\cal D}_q$ into the $\bk$-algebra $M_N(\bk)$. This implies that one has
\[
\sum^{N-1}_{k=0} D_N^{N-1-k}H^{N-1}_N D^k_N=\sum^{N-1}_{r=0}a_r D^r_N H^r_N.
\]
On the other hand one easily obtains by computation
\[
\sum^{N-1}_{k=0}D^{N-1-k}_N H^{N-1}_N D^k_N=[N-1]_q!\bbbone
\]
It follows that $(a_0-[N-1]_q!)\bbbone + \sum^{N-1}_{r=1} a_r D^r_NH^r_N=0$
which implies $a_0=[N-1]_q!$ and $a_r=0$ for $1\leq r\leq N-1$ since
$D^r_NH^r_N$ is a diagonal matrix with the first $(N-r)$ diagonal elements
invertible (by Assumption $(A_1)$) and the last $r$ diagonal elements vanish.
$\Box$\\

Let us now say a few words on the case where $\bk$ is a field. In this case
a\linebreak[4] $N$-differential module $E$ will be called a $N$-{\sl
differential vector space} and the $H_{(n)}(E)$ are vector spaces.

\begin{lemma}
If $E$ is a finite-dimensional $N$-differential vector space then, for each
integer $n$ with
$1\leq n\leq N-1$, the vector spaces $H_{(n)}(E)$ and $H_{(N-n)}(E)$ are
isomorphic.
\end{lemma}

\noindent {\bf Proof.} One has the isomorphisms of vector spaces\\
$E\simeq \ker(d^n)\oplus \im (d^n)=Z_{(n)}(E) \oplus B_{(N-n)}(E)$ and\\
$E\simeq \ker (d^{N-n})\oplus \im (d^{N-n})=Z_{(N-n)}(E) \oplus B_{(n)}(E)$.\\
On the other hand, one also has the isomorphisms of vector spaces
$Z_{(n)}(E)\simeq H_{(n)}(E)\oplus B_{(n)}(E)$ and $Z_{(N-n)}(E)\simeq
H_{(N-n)}(E)\oplus B_{(N-n)}(E)$. The isomorphism $H_{(n)}(E)\simeq
H_{(N-n)}(E)$ follows. $\Box$\\

\noindent {\bf Warning:} One must be aware of the fact that the above
isomorphisms are not canonical; in particular, they do not allow to replace the
exact hexagons of Lemma 1 and Lemma 2 by exact triangles.\\

In the case where $(E,d)$ is a finite dimensional $N$-differential vector space
over $\bk=\mathbb R$ or $\mathbb C$ one shows easily (see e.g. in \cite{Gre})
by decomposing $E$ into irreductible factors for $d$ that one has an
isomorphism $E\simeq \oplus^N_{n=1}{\mathbf k}^n\otimes {\mathbf k}^{m_n}$,
$d\simeq \oplus^N_{n=2} D_n\otimes Id_{{\mathbf k}^{m_n}}$ with

\[
D_n=\left(
\begin{array}{ccccccc}
0 & 1 & 0 & . & . & . & 0\\
. & . & . & . &  & &.\\
. &   & . & . & . &  & .\\
. &   &   & . & . & . & .\\
. &  &   &   & .  & . & 0\\
. & & & & &.& 1\\
0 & . &. &. & . & . & 0\\
\end{array}
\right)\in M_n(\bk)
\]
where the {\sl multiplicities}\  $m_n$, $n\in \{1,\dots,N\}$, are invariants of
$(E,d)$ with\linebreak[4] $\sum^N_{n=1}n m_n=\dim (E)$. Notice that one has
$m_N\geq 1$ whenever $d^{N-1}\not=0$. Notice also that the above decomposition
of $d$ is its {\sl Jordan normal-form}. In terms of the multiplicities, one can
easily compute the dimensions of the vector spaces $H_{(k)}(E)$. The result is
given by the following proposition.

\begin{proposition}
Let $E$ be a finite dimensional $N$-differential vector space over $\mathbb R$
or $\mathbb C$ with multiplicities $m_n$, $n\in \{1,2,\dots,N\}$, then one has
for each integer $k$ with $1\leq k\leq N/2$
\[
\dim H_{(k)}(E)=\dim H_{(N-k)}(E)=\sum^{k-1}_{r=1}
r(m_r+m_{N-r})+k\sum^{N-k}_{s=k}m_s,
\]
the first summation vanishing for $k=1$.
\end{proposition}
Although easy, that kind of results is useful for applications e.g. in
conformal field theory \cite{D-VT}.
\newpage

\section{$N$-complexes and (co)simplicial modules}

In this section we study the $\mathbb Z$-graded case. A $N$-{\sl complex of
modules} \cite{Kapr} is a $N$-differential module which is a $\mathbb Z$-graded
module $E=\oplusinf_{n\in \mathbb Z}E^n$ with an homogeneous $N$-differential
$d$ of degree 1 or $-1$. When $d$ is of degree 1, $E$ is refered to as {\sl a
cochain $N$-complex} and when $d$ is of degree $-1$, $E$ is refered to as a
{\sl chain $N$-complex}. One passes easily from the cochain case to the chain
case. Here we adopt the cochain language and therefore in the following a
$N$-{\sl complex}, without other specification, always means a cochain
$N$-complex of modules. An ordinary complex, or simply a complex, will mean an
ordinary cochain complex of modules that is it corresponds to the case where
$N=2$. If $E$ is a $N$-complex then the $H_{(m)}(E)$ are $\mathbb Z$-graded
modules;
$H_{(m)}(E) = \oplusinf_{n\in \mathbb Z}H^n_{(m)}(E)$ with $H_{(m)}^n(E) = \ker
(d^m:E^n\rightarrow E^{n+m})/d^{N-m}(E^{n+m-N})$. In this case the hexagon
$(\calh^{\ell,m})$ of Lemma 1 splits into long exact sequences
$(\cals^{\ell,m}_p),\  p\in \mathbb Z$\\

\[
\begin{array}{lllll}
& & & &\dots \rightarrow H^{Nr+p}_{(m)}(E)\stackrel{[i]^\ell}{\hbox to
12mm{\rightarrowfill}} H^{Nr+p}_{(\ell+m)}(E) \stackrel{[d]^m}{\hbox to
12mm{\rightarrowfill}}H^{Nr+p+m}_{(\ell)}(E)\\
\\
(\cals^{\ell,m}_p) & && & \stackrel{[i]^{N-(\ell+m)}}{\hbox to
12mm{\rightarrowfill}}H^{Nr+p+m}_{(N-m)}(E)\stackrel{[d]^\ell}{\hbox to
12mm{\rightarrowfill}}H^{Nr+p+\ell+m}_{(N-(\ell+m))}(E)\\
\\
& && &\stackrel{[i]^m}{\hbox to
12mm{\rightarrowfill}}H^{Nr+p+\ell+m}_{(N-\ell)}(E)
\stackrel{[d]^{N-(\ell+m)}}{\hbox to
12mm{\rightarrowfill}}H^{N(r+1)+p}_{(m)}(E) \stackrel{[i]^\ell}{\hbox to
12mm{\rightarrowfill}} \dots
\end{array}
\]
\\
\noindent One has $(\cals^{\ell,m}_p)=(\cals^{\ell,m}_{p+N})$ which means that
$(\cals^{\ell,m}_p)$ does only depend on the class of $p$ in $\mathbb
Z_N=\mathbb Z/N\mathbb Z$. Thus for fixed $\ell, m$ as in Lemma 1 one has only
$N$ independent long exact sequences $(\cals^{\ell,m}_p)$ which are obtained
for instance by taking $p$ in $\{0,1,\dots,N-1\}$, however it will be
convenient in the following to keep the redundant indexation by $p\in \mathbb
Z$.\\

Let $E$ and $E'$ be $N$-complexes, an {\sl homomorphism of $N$-complexes} of
$E$ into $E'$ is an homomorphism of $N$-differential modules
$\varphi:E\rightarrow E'$ which is homogeneous of degree 0, (i.e.
$\varphi(E^n)\subset E^{\prime n}$). Such an homomorphism of $N$-complexes
induces module-homomorphisms $\varphi_\ast:H^n_{(m)}(E)\rightarrow
H^n_{(m)}(E')$ for $n\in \mathbb Z$ and $1\leq m\leq N-1$. Let $0\rightarrow E
\stackrel{\alpha}{\rightarrow} F \stackrel{\beta}{\rightarrow} G\rightarrow 0$
be a short exact sequence of $N$-complexes, then the hexagon $(\calh_n)$ of
Lemma 2 splits into  long exact sequences $(\cals_{n,p}), p\in\mathbb Z$\\

\[
\begin{array}{lllll}
& & & &\dots \rightarrow H^{Nr+p}_{(n)}(E)\stackrel{\alpha_\ast}{\hbox to
12mm{\rightarrowfill}} H^{Nr+p}_{(n)}(F) \stackrel{\beta_\ast}{\hbox to
12mm{\rightarrowfill}}H^{Nr+p}_{(n)}(G)\\
\\
(\cals_{n,p}) & && &
\stackrel{\partial}{\rightarrow}H^{Nr+p+n}_{(N-n)}(E)\stackrel{\alpha_\ast}{\hbox to 12mm{\rightarrowfill}}H^{Nr+p+n}_{(N-n)}(F)\stackrel{\beta_\ast}{\hbox to 12mm{\rightarrowfill}} H^{Nr+p+n}_{(N-n)}(G)\\
\\
& && &\stackrel{\partial}{\rightarrow}H^{N(r+1)+p}_{(n)}(E)
\stackrel{\alpha_\ast}{\hbox to 12mm{\rightarrowfill}} \dots
\end{array}
\]
\\
\noindent One has again $(\cals_{n,p})=(\cals_{n,p+N})$ so there are $N$ such
sequences for a given integer $n$ with $1\leq n\leq N-1$.\\

Let $E$ be a $N$-complex, we associate to $E$ the following ordinary complexes
$C_{m,p}$, $1\leq m\leq N-1,\ p\in \mathbb Z$\\

\[
C_{m,p}\ \ \ \ \ \ \dots \stackrel{d^{N-m}}{\gf} E^{Nr+p} \stackrel{d^m}{\gf}
E^{Nr+p+m} \stackrel{d^{N-m}}{\gf} E^{N(r+1)+p} \stackrel{d^m}{\gf} \dots
\]
\\
\noindent that is $C_{m,p}=\oplusinf_{n\in \mathbb Z} C^n_{m,p}$ with
$C^{2r+p}_{m,p}= E^{Nr+p}$ and $C^{2r+p+1}_{m,p}=E^{Nr+p+m}$, the differential
of $C_{m,p}$ being given alternatively by $d^m$ and $d^{N-m}$. Again one has
$C_{m,p}=C_{m,p+N}$ so, for fixed $m$ in $\{1,\dots,N-1\}$, there are
generically only $N$ distinct such complexes which are obtained for instance by
taking $p$ in $\{0,1,\dots,N-1\}$ but it will be convenient in the following to
keep the redundant indexation by $p\in \mathbb Z$. Notice that one also has
$C_{m,p}=C_{N-m,p+m}$ (so there may be less than $N$ distinct relevant $p$).
Let $\ell$ and $m$ be integers with $\ell\geq 1, m\geq 1, \ell+m\leq N-1$ and
consider the following commutative diagram\\

\[
\begin{array}{lllllll}
\dots \stackrel{d^{N-m}}{\gf} & E^{Nr+p} & \stackrel{d^m}{\gf} & E^{Nr+p+m} &
\stackrel{d^{N-m}}{\gf} & E^{N(r+1)+p} & \stackrel{d^m}{\gf}\dots\\
& \downarrow I &  & \downarrow d^\ell & & \downarrow I &\\
\dots \stackrel{d^{N-(\ell+m)}}{\gf} & E^{Nr+p} & \stackrel{d^{\ell+m}}{\gf} &
E^{Nr+p+\ell+m} & \stackrel{d^{N-(\ell+m)}}{\gf} & E^{N(r+1)+p} &
\stackrel{d^{\ell+m}}{\gf}\dots
\end{array}
\]
\\
\noindent where $I$ denotes the identity mappings of the corresponding modules.
The first row is the complex $C_{m,p}$ and the second row is the complex
$C_{\ell+m,p}$ so the columns define an homomorphism of complexes
$\Delta^\ell:C_{m,p}\rightarrow C_{\ell+m,p}$; this notation is justified by
the fact that one has $\Delta^\ell=(\Delta^1)^\ell,
\Delta^\ell=\Delta^r\Delta^s$ if $\ell=r+s$ with $r\geq 1, s\geq 1$. In
(co)homology, $\Delta^\ell$ induces homomorphisms
$\Delta^\ell_\ast:H^n(C_{m,p})\rightarrow H^n(C_{\ell+m,p})$;
$\Delta^\ell_\ast:H^\ast(C_{m,p})\rightarrow H^\ast(C_{\ell+m,p})$ is an
homomorphism of graded modules.

\begin{lemma}
The following conditions $(a)$, $(b)$ and $(c)$ are equivalent.\\
$(a)$ $H^{Nr+p+m}_{(\ell)}(E)=0$ and $H^{Nr+p+\ell+m}_{(N-\ell)}(E)=0, \ \
\forall r\in \mathbb Z$\\
$(b)$ $[i]^\ell:H^{Nr+p}_{(m)}(E)\rightarrow H^{Nr+p}_{(\ell+m)}(E)$ and
$[d]^\ell:H^{Nr+p+m}_{(N-m)}(E)\rightarrow  H^{Nr+p+\ell+m}_{(N-(\ell+m))}(E)$
are isomorphisms, $\forall r\in \mathbb Z$.\\
$(c)$ $\Delta^\ell_\ast:H^\ast(C_{m,p})\rightarrow H^\ast(C_{\ell+m,p})$ is an
isomorphism.
\end{lemma}
\noindent {\bf Proof.} The equivalence $(a)\Leftrightarrow (b)$ follows from
the exactness of the sequence $(\cals^{\ell,m}_p)$. The equivalence
$(b)\Leftrightarrow (c)$ follows from
$H^{Nr+p}_{(m)}(E)=H^{2r+p}(C_{m,p})$,\linebreak[4]
$H^{Nr+p}_{(\ell+m)}(E)=H^{2r+p}(C_{\ell+m,p})$ and $\Delta^\ell_\ast=[i]^\ell$
on $H^{2r+p}(C_{m,p})=H^{Nr+p}_{(m)}(E)$ and from
$H^{Nr+p+m}_{(N-m)}(E)=H^{2r+p+1}(C_{m,p}),
H^{Nr+p+\ell+m}_{(N-(\ell+m))}(E)=H^{2r+p+1}(C_{\ell+m,p})$ and
$\Delta^\ell_\ast=[d]^\ell$ on $H^{2r+p+1}(C_{m,p})=H^{Nr+p+m}_{(N-m)}(E)$.
$\Box$

\begin{corollary}
The following conditions $(i)$, $(ii)$ and $(iii)$ are equivalent.\\
$(i)$ $H_{(m)}^n(E)=0$ if $n\not =p$ $(mod N)$ and $n+m\not= p$     $(mod
N)$,\\
$(ii)$ $\Delta_\ast:H^\ast(C_{m,p})\rightarrow H^\ast(C_{m+1,p})$ is an
isomorphism, $\forall m\in \{1,\dots,N-2\}$,\\
$(iii)$ $H^{Nr+p+m}_{(1)}(E)=0$ and $H^{Nr+p+m+1}_{(N-1)}(E)=0$, $\forall m\in
\{1,\dots,N-2\}$ and $\forall r\in \mathbb N$.
\end{corollary}

In the following, we shall be mainly concerned with positive $N$-complexes,
that is with the $\mathbb N$-graded case $E=\oplus_{n\in \mathbb N} E^n$ (or
$E^n=0$ if $n<0$). There is a standard way to produce positive complexes
starting from (co)simplicial modules, (see e.g. \cite{J-LL}, \cite{Weib}). A
{\sl pre-cosimplicial module} (or {\sl semi-cosimplicial} in the terminology of
\cite{Weib}) is a sequence of modules $(E^n)_{n\in \mathbb N}$ together with
{\sl coface homomorphisms} ${\mathfrak f}_i: E^n\rightarrow E^{n+1},\ i\in
\{0,1,\dots,n+1\}$, satisfying\\
$({\Fg})$\hspace{2cm} ${\fg}_j {\fg}_i={\fg}_i{\fg}_{j-1}$\hspace{.5cm} if
$i<j$.\\
Given a pre-cosimplicial module $(E^n)$, one associates to it a positive
complex $(E,d)$ by setting $E=\oplus_{n\in \mathbb N} E^n$ and $
d=\sum^{n+1}_{i=0} (-1)^i \fg_i: E^n\rightarrow E^{n+1}$.\\
One verifies that $d^2=0$ is implied by the coface relations $(\Fg)$.\\
In the rest of this section $(E^n)$ is a fixed pre-cosimplicial module, $d$
will always denote the above differential and $q$ is a fixed element of $\bk$
such that $[N]_q=0$, (i.e. $\bk$ and $q\in \bk$ satisfy assumption $(A_0)$ of
\S 2). Furthermore, the $H^n=H^n(E)=\ker(d:E^n\rightarrow E^{n+1})/dE^{n-1}$
will denote the (co)homology of $(E,d)$ and will be refered to as {\sl the
cohomology of the pre-cosimplicial module} $(E^n)$.\\
Let us define two sequences $d_m,\delta_m$ ($m\in \mathbb N$) of homogeneous
endomorphisms of degree 1 of $E$ by setting\\
$\delta_m=\sum^{n-m+1}_{i=0} q^i \fg_i:E^n\rightarrow E^{n+1}$ for $n\geq m-1$
and $\delta_m(E^n)=0$ for $n<m-1$\\
$d_m=\delta_{m+1}+q^{n-m+1}\left(\sum^m_{r=0}(-1)^r
\fg_{n-m+r+1}\right):E^n\rightarrow E^{n+1}$ for $n\geq m-1$\\
$d_m=d:E^n\rightarrow E^{n+1}$ for $n<m-1$ and $d_0=\delta_0$.\\
Notice that $d_{n+1}=d$ on $E^n$ and that therefore $d_m=d$ on $E^n$ for $n\leq
m-1$. For $n\geq m-1\geq 0$ one can also write\\
$d_m=\delta_m-q^{n-m+1}(\sum^{m-1}_{r=0}(-1)^r \fg_{n-m+r+2})$ :
$E^n\rightarrow E^{n+1}$.

\begin{lemma}
Let $n$ and $m$ be elements of $\mathbb N$ with $n\geq m$. Then, on $E^n$ one
has the following identities:\\
$(i)$ $d^{p+1}_m=\delta^p_{m+1} (\delta_{m+1} +[p+1]_q
q^{n-m+1}\sum^m_{r=0}(-1)^r\fg_{n-m+r+1}),\ \ \forall p\in \mathbb N$\\
$(ii)$ $\delta^{p+1}_m=\delta^p_{m+1} (\delta_{m+1} +[p+1]_q q^{n-m+1}
\fg_{n-m+1}), \ \ \forall p\in \mathbb N.$\\
Furthermore the following identity holds on $E$ for any $m\in \mathbb N$\\
$(iii)$ $d_{m+1} d_m=\delta^2_{m+1}$.
\end{lemma}

\noindent{\bf Proof.} By induction on $p$. From the definitions it follows that
$(i)$ and $(ii)$ are true for $p=0$. The passage  from $p$ to $p+1$ follows
from $(\Fg)$ which implies
\[
\begin{array}{ll}
&(\delta_{m+1}+[p+1]_qq^{n-m+2}\sum^m_{r=0}(-1)^r\fg_{n-m+r+2})d_m =\\
\\
& \delta_{m+1}(\delta_{m+1}+[p+2]_qq^{n-m+1}\sum^m_{r=0}(-1)^r\fg_{n-m+r+1})
\end{array}
\]
and
\[
(\delta_{m+1}+[p+1]_qq^{n-m+2}\fg_{n-m+2})\delta_m=\delta_{m+1}(\delta_{m+1}+[p+2]_q q^{n-m+1}\fg_{n-m+1}).
\]
For $(iii)$ one first notices that on $E^n$ with $n\leq m-1$ one has
$d_{m+1}d_m=d^2=0$ and $\delta_{m+1}=0$. On the other hand if $n\geq m$ one has
on $E^n$
\[
\begin{array}{l}
d_{m+1}d_m =\\
\\
= (\delta_{m+1}-q^{n-m+1}\sum^m_{r=0}(-1)^r\fg_{n-m+r+2})(\delta_{m_+1}
+ q^{n-m+1}\sum^m_{s=0}(-1)^s\fg_{n-m+s+1})\\
\\
=
\delta^2_{m+1}+q^{n-m+1}(\delta_{m+1}-\delta_{m+2})\sum^m_{r=0}(-1)^r\fg_{n-m+r+1}\\
\\
-
(q^{n-m+1})^2(\sum^m_{r=0}(-1)^r\fg_{n-m+r+2})(\sum^m_{s=0}(-1)^s\fg_{n-m+s+1})\\
\\
= \delta^2_{m+1} +
(q^{m-m+1})^2(\sum^{m+1}_{r=0}(-1)^r\fg_{n-m+r+1})(\sum^m_{s=0}(-1)^s\fg_{n-m+s+1})\\
\\
= \delta^2_{m+1}. \Box
\end{array}
\]
It is worth noticing that the assumption $(A_0)$ (i.e. $[N]_q=0$) is not needed
for this lemma. It is needed however for next corollary.

\begin{corollary}
One has the following identities on $\oplus_{n\geq m} E^n$\\
$(a)$ $d^{N-1}_m=\delta^{N-2}_{m+1} d_{m+1}$\\
$(b)$ $d_{m-1}\delta^{N-2}_m =\delta^{N-2}_{m+1}d_{m+1}\ \ \mbox{for}\ \ m\geq
1$.
\end{corollary}

\noindent {\bf Proof.} Applying the identity $(i)$ of Lemma 8 on $E^n$ for
$p=N-2$, one obtains
\[
\begin{array}{ll} &
d^{N-1}_m=\delta^{N-2}_{m+1}(\delta_{m+1}+[N-1]_qq^{n-m+1}\sum^m_{r=0}(-1)^r\fg_{n-m+r+1})=\\
\\
&
\delta^{N-2}_{m+1}(\delta_{m+1}-q^{n-m}\sum^m_{r=0}(-1)^r\fg_{n-m+r+1})=\delta^{N-2}_{m+1}d_{m+1}
\end{array}
\]
by using $0=[N]_q=1+q[N-1]_q$. This proves the identity $(a)$. Similarly, one
proves $(b)$ by using the identity $(ii)$ of Lemma 8 on $E^n$ for $p=N-2$ and
the fact that, by $(\Fg)$, one has $\fg_{n+N-2-m+r+2}\ \delta^{N-2}_m =
\delta^{N-2}_{m+1} \fg_{n-m+r+2}$ on $E^n$, (and of course again $[N]_q=0$, or
equivalently $q[N-1]_q=-1)$. $\Box$

\begin{corollary}
The operators $d_m$ and $\delta_m$ satisfy $d^N_m=0$ and $\delta^N_m=0$ for any
$m\in \mathbb N$.
\end{corollary}

\noindent{\bf Proof.} For $\delta_m$, this is classical (see e.g. \cite{Kapr}),
namely one shows by induction that ($\Fg$) implies
$\delta^p_m=[p]_q!\sum_{i_1>\dots>i_p} q^{i_1+\dots +i_p-\frac{p(p-1)}{2}}
\fg_{i_1}\dots \fg_{i_p}$ where the limit of the sum depends on $m$ and on the
degree of the element on which one applies the equality. Thus $[N]_q=0$ implies
$\delta^N_m=0$.\\
For $d_m$, one obtains by applying $d_m$ on the right to both sides of the
identity $(a)$ of Corollary 3 and by using the identity $(iii)$ of Lemma 8:
$d^N_m=\delta^{N-2}_{m+1} d_{m+1}d_m=\delta^N_{m+1}=0$ on $\oplus_{n\geq
m-1}E^n$. On the other hand if $n<m-1$ then on $E^n$ one has $d^2_m=d^2=0$
which achieves the proof of $d^N_m=0$ since $N\geq 2$. $\Box$
\\

Thus, for each $p\in \mathbb N$, ($E,d_p$) is a $N$-complex. It follows from
Corollary 3 that one has on $E^n$,  $d^{N-1}_p=\delta^{N-2}_{p+1} d_{p+1}$ for
$n\geq p\geq 0$ and $d^{N-1}_{p+1}=d_p\delta^{N-2}_{p+1}$ for $n\geq p+1\geq
1$. In other words one has the following commutative diagram $(\Psi)$\\

\[
(\Psi)\ \
\begin{diagram}
\node{E^n}\arrow{e,t}{d_{p+1}}\arrow{s,r}{I}
\node{E^{n+1}}\arrow{e,t}{d^{N-1}_{p+1}}\arrow{s,r}{\delta^{N-2}_{p+1}}
\node{E^{n+N}}\arrow{e,t}{d_{p+1}}\arrow{s,r}{I}
\node{E^{n+N+1}}\arrow{e,t}{d^{N-1}_{p+1}}\arrow{s,r}{\delta^{N-2}_{p+1}}
\node{\dots}\\
\node{E^n} \arrow{e,t}{d^{N-1}_p}
\node{E^{n+N-1}}\arrow{e,t}{d_p}
\node{E^{n+N}}\arrow{e,t}{d^{N-1}_p}
\node{E^{n+2N-1}}\arrow{e,t}{d_p}
\node{\dots}
\end{diagram}
\]
\\
\noindent for $n\geq p\geq 0$.

\begin{corollary}
The above commutative diagram induces homomorphisms
\[
\begin{array}{ll}
\Psi_\ast:H^{n+1}_{(N-1)}(E,d_{p+1}) & \rightarrow H^{n+N-1}_{(1)}(E,d_p),\\
\\
\Psi_\ast:H^{n+N}_{(1)}(E,d_{p+1}) & \rightarrow H^{n+N}_{(N-1)}(E,d_p)
\end{array}
\]
for any integer $n,p$ with $n\geq p\geq 0$ and, for $p=0$
\[
\Psi_\ast:H^0_{(1)}(E,d_1)\rightarrow H^0_{(N-1)}(E,d_0)
\]
\end{corollary}
Notice that the last homomorphism (``for $p=0$") follows from the obvious fact
that if $p=0$ or if $p=1$, one has $H^0_{(m)}(E,d_p)=Z^0_{(m)}(E,d_p)$.\\

In order to avoid unessential complications in low degrees, we shall truncate
the $N$-complex $(E,d_p)$ by setting $E_p=\oplus_{n\geq p-1}E^n$ with the
convention $E^{-1}=0$ for $p=0$, (thus $E_0=E_1=E$). Clearly $(E_p,d_p)$ is
still a $N$-complex and our aim in the remaining part of this section is to
rely its generalized cohomology to the cohomology of the pre-cosimplicial
module ($E^n$) or, more precisely to the cohomology of the truncated complex
($E_p,d$); one obviously has $H^n(E)=H^n(E_p,d)$, $\forall n\geq p$.

\begin{theorem}
The following diagram of homomorphisms
\[
\begin{diagram}
\node{\dots}\arrow{e,t}{d}
\node{E^{2r+p-1}}\arrow{e,t}{d}\arrow{s,r}{{\delta^{N-2}_{2+p-1}}
\dots \delta^{N-2}_{2r+p-1}}
\node{E^{2r+p}}\arrow{e,t}{d}
\arrow{s,r}{{\delta^{N-2}_{2+p}}\dots \delta^{N-2}_{2r+p}}
\node{E^{2(r+1)+p-1}}\arrow{e,t}{d}\arrow{s,r}{{\delta^{N-2}_{2+p-1}}\dots
\delta^{N-2}_{2(r+1)+p-1}}\node{\dots}\\
\node{\dots} \arrow{e,t}{d^{N-1}_p}\node{E^{Nr+p-1}}\arrow{e,t}{d_p}
\node{E^{Nr+p}}\arrow{e,t}{d^{N-1}_p}
\node{E^{N(r+1)+p-1}}\arrow{e,t}{d_p}\node{\dots}
\end{diagram}
\]
starting with
\[
\begin{diagram}
\node{E^{p-1}}\arrow{s,r}{I}\arrow{e,t}{d}\node{E^p}\arrow{s,r}{I}\arrow{e,t}{d} \node{E^{p+1}}\arrow{s,r}{\delta^{N-2}_{p+1}}\arrow{e,r}{d}\node{\dots}\\
\node{E^{p-1}}\arrow{e,t}{d_p}\node{E^p}\arrow{e,t}{d^{N-1}_p}\node{E^{N+p-1}}\arrow{e,t}{d_p}\node{\dots}
\end{diagram}
\]
is commutative, $\forall p\in \mathbb N$.
\end{theorem}

\noindent{\bf Proof.} By Corollary 3 (b), one has $d_p
\delta^{N-2}_{2+p-1}\dots \delta^{N-2}_{2r+p-1}=\delta^{N-2}_{2+p}\dots
\delta^{N-2}_{2r+p}d_{2r+p}$ and by Corollary 3 (a) and (b), one has
\[
d^{N-1}_p \delta^{N-2}_{2+p} \dots
\delta^{N-2}_{2r+p}=\delta^{N-2}_{2+p-1}\dots\delta^{N-2}_{2(r+1)+p-1}d_{2r+p+1}.
\]
The result follows then from the fact that, for any $n\in \mathbb N$,
$d_{n+1}=d$ on $E^n$. $\Box$\\

In Theorem 1, the first row of the diagram is the truncated complex $(E_p,d)$
whereas the second row is the complex $C_{1,p-1}(E_p,d_p)$ associated with
the\linebreak[4] $N$-complex $(E_p,d_p)$ as defined just before Lemma 7 so the
columns of the diagram define an homomorphism of complexes
$\Phi_{1,p}:(E_p,d)\rightarrow C_{1,p-1}(E_p,d_p)$. By combining this
homomorphism with the homomorphism $\Delta$ defined before Lemma 7, one obtains
for each integer $m$ with $1\leq m\leq N-1$ an homomorphism of complexes
$\Phi_{m,p}=\Delta^{m-1}\circ \Phi_{1,p}:(E_p,d)\rightarrow
C_{m,p-1}(E_p,d_p)$. The situation is thus summarized by the following
commutative diagram\\

\[
\begin{diagram}
\node{}\arrow{e,t}{d}
\node{E^{2r+p-1}}\arrow{s,r}{\delta^{N-2}_{2+p-1}\dots\delta^{N-2}_{2r+p-1}}\arrow{e,t}{d}
\node{E^{2r+p}}\arrow{s,r}{\delta^{N-2}_{2+p}\dots\delta^{N-2}_{2r+p}}\arrow{e,t}{d}
\node{E^{2(r+1)+p-1}}\arrow{s,r}{\delta^{N-2}_{2+p-1}\dots
\delta^{N-2}_{2(r+1)+p-1}}\arrow{e,t}{d}
\node{\dots}\\
\node{}\arrow{e,t}{d^{N-1}_p}
\node{E^{Nr+p-1}}\arrow{s,r}{I}\arrow{e,t}{d_p}
\node{E^{Nr+p}}\arrow{s,r}{d^{m-1}_p}\arrow{e,t}{d^{N-1}_p}
\node{E^{N(r+1)+p-1}}\arrow{s,r}{I}\arrow{e,t}{d_p}
\node{\dots}
\\
\node{}\arrow{e,t}{d^{N-m}_p}
\node{E^{Nr+p-1}}\arrow{e,t}{d^m_p}
\node{E^{Nr+m+p-1}} \arrow{e,t}{d^{N-m}_p}
\node{E^{N(r+1)+p-1}}\arrow{e,t}{d^m_p}
\node{\dots}
\end{diagram}
\]
\\
\noindent starting with\\

\[
\begin{diagram}
\node{E^{p-1}}\arrow{s,r}{I}\arrow{e,t}{d}
\node{E^p}\arrow{s,r}{I}\arrow{e,t}{d}
\node{E^{p+1}}\arrow{s,r}{\delta^{N-2}_{p+1}}\arrow{e,t}{d}
\node{\dots}\\
\node{E^{p-1}}\arrow{s,r}{I}\arrow{e,t}{d_p}
\node{E^p}\arrow{s,r}{d^{m-1}_p}\arrow{e,t}{d^{N-1}_p}
\node{E^{N+p-1}}\arrow{s,r}{I}\arrow{e,t}{d_p}
\node{\dots}\\
\node{E^{p-1}}\arrow{e,t}{d^m_p}
\node{E^{m+p-1}}\arrow{e,t}{d^{N-m}_p}
\node{E^{N+p-1}}\arrow{e,t}{d^m_p}
\node{\dots}
\end{diagram}
\]
\\
\begin{corollary}
The homomorphism $\Phi_{m,p}$ induces canonical homomorphisms of modules
\[
\begin{array}{lll}
\Phi_{\ast m,p}:H^{2r+p-1}(E) & \longrightarrow H^{Nr+p-1}_{(m)}(E_p,d_p) &
$\mbox{for}$\ \ r\in \mathbb N\ $\mbox{with}$\ \ r\geq 1\\
\\
\Phi_{\ast m,p}:H^{2r+p}(E) & \longrightarrow H^{Nr+m+p-1}_{(N-m)}(E_p,d_p) &
$\mbox{for}$\ \ r\in \mathbb N,
\end{array}
\]
for any $p\in \mathbb N$ and $m\in \{1,\dots,N-1\}$. Furthermore, for $p=1$ one
also has a canonical homomorphism $\Phi_{\ast m,1}:H^0(E)\rightarrow
H^0_{(m)}(E,d_1)$, (i.e. one can take $r=0$), for any $m\in \{ 1,\dots,N-1\}$.
\end{corollary}

\noindent{\bf Proof.} This is just the homomorphisms in (co)homology associated
to the homomorphism of complexes. $\Box$\\

Before adding some more assumptions, let us make contact between Theorem 1 and
the diagram $(\Psi)$. Since $d_{p+1}=d$ on $E^p$, in the case $n=p$ one can
extend on the left the diagram $(\Psi)$ by adding the beginning of the complex
$(E,d)$. More precisely let us fill up the first quadrant, with the following
commutative diagram \\

\[
\begin{diagram}
\node{}\arrow{e,!}
\node{}\arrow{e,!}\arrow{s,r}{I}
\node{}\arrow{e,!}\arrow{s,r}{I}
\node{}\arrow{e,!}\arrow{s,r}{I}
\node{}\arrow{e,!}\arrow{s,r}{I}
\\
\node{}\arrow{e,t}{d}
\node{E^{p-2}}\arrow{e,t}{d}\arrow{s,r}{I}
\node{E^{p-1}}\arrow{e,t}{d}\arrow{s,r}{I}
\node{E^p}\arrow{e,t}{d=d_{p+1}}\arrow{s,r}{I}
\node{E^{p+1}}\arrow{e,t}{d^{N-1}_{p+1}}\arrow{s,r}{\delta^{N-2}_{p+1}}
\node{}
\\
\node{}\arrow{e,t}{d}
\node{E^{p-2}}\arrow{e,t}{d}\arrow{s,r}{I}
\node{E^{p-1}}\arrow{e,t}{d=d_p}\arrow{s,r}{I}
\node{E^p}\arrow{e,t}{d^{N-1}_p}\arrow{s,r}{\delta^{N-2}_p}
\node{E^{N+p-1}}\arrow{e,t}{d_p}\arrow{s,r}{I}
\node{}
\\
\node{}\arrow{e,t}{d}
\node{E^{p-2}}\arrow{e,t}{d=d_{p-1}}\arrow{s,r}{I}
\node{E^{p-1}}\arrow{e,t}{d^{N-1}_{p-1}}\arrow{s,r}{\delta^{N-2}_{p-1}}
\node{E^{N+p-2}}\arrow{e,t}{d_{p-1}}\arrow{s,r}{I}
\node{E^{N+p-1}}\arrow{e,t}{d^{N-1}_{p-1}}\arrow{s,r}{\delta^{N-2}_{p-1}}
\node{}
\\
\node{}\arrow{e,t}{d=d_{p-2}}
\node{E^{p-2}}\arrow{e,t}{d^{N-1}_{p-2}}\arrow{s,r}{\delta^{N-2}_{p-2}}
\node{E^{N+p-3}}\arrow{e,t}{d_{p-2}}\arrow{s,r}{I}
\node{E^{N+p-2}}\arrow{e,t}{d^{N-1}_{p-2}}\arrow{s,r}{\delta^{N-2}_{p-2}}
\node{E^{2N+p-3}}\arrow{e,t}{d_{p-2}}\arrow{s,r}{I}
\node{}
\end{diagram}
\]
\\

\vspace{2cm}

with corner\\

\[
\begin{diagram}
\node{\vdots}\arrow{e,!}\arrow{s,r}{I}
\node{\vdots}\arrow{e,!}\arrow{s,r}{I}
\node{\vdots}\arrow{e,!}\arrow{s,r}{I}
\node{\vdots}\arrow{e,!}\arrow{s,r}{I}
\node{}\\
\node{E^0}\arrow{e,t}{d}\arrow{s,r}{I}
\node{E^1}\arrow{e,t}{d}\arrow{s,r}{I}
\node{E^2}\arrow{e,t}{d=d_3}\arrow{s,r}{I}
\node{E^3}\arrow{e,t}{d^{N-1}_3}\arrow{s,r}{\delta^{N-2}_3}
\node{\dots}\\
\node{E^0}\arrow{e,t}{d}\arrow{s,r}{I}
\node{E^1}\arrow{e,t}{d=d_2}\arrow{s,r}{I}
\node{E^2}\arrow{e,t}{d^{N-1}_2}\arrow{s,r}{\delta^{N-2}_2}
\node{E^{N+1}}\arrow{e,t}{d_2}\arrow{s,r}{I}
\node{\dots}\\
\node{E^0}\arrow{e,t}{d=d_1}\arrow{s,r}{I}
\node{E^1}\arrow{e,t}{d^{N-1}_1}\arrow{s,r}{\delta^{N-2}_1}
\node{E^N}\arrow{e,t}{d_1}\arrow{s,r}{I}
\node{E^{N+1}}\arrow{e,t}{d^{N-1}_1}\arrow{s,r}{\delta^{N-2}_1}
\node{\dots}\\
\node{E^0}\arrow{e,t}{d^{N-1}_0}
\node{E^{N-1}}\arrow{e,t}{d_0}
\node{E^N}\arrow{e,t}{d^{N-1}_0}
\node{E^{2N-1}}\arrow{e,t}{d_0}
\node{\dots}
\end{diagram}
\]
\\

\noindent The $p$-th row is the complex $C_{1,p-1}(E_p,d_p)$ on the right of
$E^{p-1}$ for $p\geq 1$ and is the beginning of the complex $(E,d)$ on the left
of $E^p$; it is a complex that we denote by $\bar C_p$. The columns between
$\bar C_{p+1}$ and $\bar C_p$ define an homomorphism of complexes $\bar
\Psi_p:\bar C_{p+1}\rightarrow \bar C_p$ which extends the commutative diagram
$(\Psi)$ when $n=p$. It is clear that the columns of the diagram of Theorem 1
are appropriate products of the $\bar\Psi_p$'s at appropriate places. In
particular, the homomorphisms $\Phi_{\ast 1,p}$ in Corollary 6 are appropriate
products of the homomorphisms $\bar\Psi_{\ast p}$ induced in (co)homology by
the $\Psi_p$:
\[
\begin{array}{ll}
\Phi_{\ast 1,p}=\bar\Psi_{\ast p}\circ\dots\circ\bar\Psi_{\ast 2r+p-1} :
H^{2r+p-1}(E) &\rightarrow H^{Nr+p-1}_{(1)}(E_p,d_p)\\
\\
\Phi_{\ast 1,p}=\bar\Psi_{\ast p}\circ\dots\circ\bar\Psi_{\ast 2r+p} :
H^{2r+p}(E) &\rightarrow H^{Nr+p}_{(N-1)}(E_p,d_p)
\end{array}
\]
Similarly when $n=p+\ell$ with $1\leq\ell\leq N-2$ one can extend on the left
the diagram $(\Psi)$ as\\

\[
\begin{diagram}
\node{\dots}\arrow{e}
\node{0}\arrow{e}\arrow{s}
\node{0}\arrow{e}\arrow{s}
\node{E^{p+\ell}}\arrow{e,t}{d_{p+1}}\arrow{s,r}{I}
\node{E^{p+\ell+1}}\arrow{e,t}{d^{N-1}_{p+1}}\arrow{s,r}{\delta^{N-2}_{p+1}}
\node{\dots}\\
\node{\dots}\arrow{e}
\node{0}\arrow{e}
\node{E^{p-1+\ell}}\arrow{e,b}{d_p}
\node{E^{p+\ell}}\arrow{e,b}{d^{N-1}_p}
\node{E^{N+p-1+\ell}}\arrow{e,b}{d_p}
\node{\dots}
\end{diagram}
\]
\\
\noindent i.e. the first row is the complex $C_{1,p+\ell}$ ($E_{p+1},d_{p+1}$)
on the right of $E^{p+\ell}$ extended by the beginning of the $0$ complex on
the left, etc. This diagram is still commutative, the rows are complexes and
the columns define an homomorphism of complexes $\bar\Psi_{p,\ell}$ which
induces in homology homorphisms\\
$\bar\Psi_{\ast p,\ell}:H^{Nr+p+\ell}_{(1)}(E_{p+1},d_{p+1})\rightarrow
H^{Nr+p+\ell}_{(N-1)}(E_p,d_p)$ for $r\geq 1$,\\
$\bar\Psi_{\ast p,\ell}: H^{Nr+p+\ell+1}_{(N-1)}(E_{p+1},d_{p+1})\rightarrow
H^{N(r+1)+p-1+\ell}_{(1)} (E_p,d_p)$ for $r\geq 0$,\\
$\bar\Psi_{\ast p,\ell}:Z^{p+\ell}_{(1)}(E_{p+1},d_{p+1})\rightarrow
H^{p+\ell}_{(N-1)}(E_p,d_p)$ and\\
$\bar\Psi_{\ast p,\ell}:0\rightarrow Z^{p-1+\ell}_{(1)}(E_p,d_p)$\\
for any $p\in \mathbb N$ and $\ell\in\{ 1,\dots,N-2\}$. This implies in
particular that one has the homomorphisms:\\
$\bar\Psi_{\ast p,\ell}\circ \dots\circ \bar\Psi_{\ast 2r+p,\ell}:0\rightarrow
H^{Nr+p-1+\ell}_{(1)}(E_p,d_p)$\\
$\bar\Psi_{\ast p,\ell}\circ \dots\circ \bar\Psi_{\ast
2r+p+1,\ell}:0\rightarrow H^{Nr+p+\ell}_{(N-1)}(E_p,d_p)$\\
for any $p\in \mathbb N$ and $\ell\in\{ 1,\dots,N-2\}$.\\

The above results are more or less optimal in the framework of pre-cosimplicial
modules under Assumption $(A_0)$ for $\bk$ and $q\in \bk$. To go further we
shall need Assumption $(A_1)$ for $\bk$ and $q\in \bk$ and we shall restrict
attention to cosimplicial modules \cite{J-LL}, \cite{Weib}. A {\sl cosimplicial
module} is a pre-cosimplicial module $(E^n)$ with coface homomorphisms $\fg_i$
as before together with {\sl codegeneracy homomorphisms}
$\sg_i:E^{n+1}\rightarrow E^n,\ i\in \{0,\dots,n\}$, satisfying\\
$(\Sg)\hspace{1.5cm} \sg_j\sg_i=\sg_i \sg_{j+1}\hspace{1cm}$ if $i\leq j$\\
and\\
$(\Sg\Fg)\hspace{1.5cm} \sg_j\fg_i=\left\{
\begin{array}{lll}
\fg_i\sg_{j-1} & \mbox{if} & i<j\\
I & \mbox{if} & i=j\  \mbox{or}\  i=j+1\\
\fg_{i-1}\sg_j & \mbox{if} & i>j+1
\end{array}
\right.
$
\\

Our aim is to prove the following theorem.
\begin{theorem}
If $(E^n)$ is a cosimplicial module and if $\bk$ and $q\in \bk$ satisfy
Assumption $(A_1)$, then the homomorphisms $\Phi_{\ast m,p}$ of Corollary 6 are
isomorphisms.
\end{theorem}
Therefore in the remaining part of this section, $(E^n)$ is a fixed
cosimplicial module, $q$ is a fixed element of $\bk$ such that $[N]_q=0$ and
$[n]_q$ is invertible for $1\leq n\leq N-1$ and we use the notations previously
introduced.\\
It follows from Corollary 4 that $\delta^N_m=0$. For $m=0$, one has
$\delta_0=d_0$ by definition. On the other hand $\delta_{p+1}$ vanishes on
$E^n$ for $n<p$ therefore we consider the action of $\delta_{p+1}$ on
$E_{p+1}=\oplus_{n\geq p} E^n$.

\begin{lemma}
For any $p\in \mathbb N$ and for any integer $m$ with $1\leq m\leq N-1$, one
has $H_{(m)}(E_{p+1},\delta_{p+1})=0$.
\end{lemma}
\noindent {\bf Proof.} Define the endomorphism $\chi_{p+1}$ of $E_{p+1}$ by
$\chi_{p+1}(E^p)=0$ and by $\chi_{p+1}=q^{-n+p}\sg_{n-p}:E^{n+1}\rightarrow
E^n$ for $n\geq p$. By using the relations $(\Sg \Fg)$ one obtains
$\chi_{p+1}\delta_{p+1}-q^{-1}\delta_{p+1}\chi_{p+1}=Id_{E_{p+1}}$ which
implies the result in view of Lemma 5 and of the fact that $\bk$ and $q^{-1}\in
\mathbb \bk$ satisfy Assumption $(A_1)$. $\Box$

\begin{lemma}
Let $p,n$ and $m$ be integers with $n\geq p\geq 0$ and $1\leq m\leq N-1$ and
let $x$ be an element of $E^n$. Then the following conditions are equivalent:\\
$(i)$ $d^m_p(x)=0$\\
$(ii)$ there exists $y\in E^{n+m-N}$ such that
\[
\delta_{p+1}(x-d^{N-m}_py)+[m]_qq^{n-p+1}\sum^p_{i=0}\fg_{n-p+i+1}(x-d^{N-m}_py)=0
\]
\end{lemma}

\noindent{\bf Proof.} By using the identity $(i)$ of Lemma 8, $d^m_p(x)=0$ is
equivalent to
\[
\delta^{m-1}_{p+1}(\delta_{p+1}(x)+[m]_qq^{n-p+1}\sum^p_{i=0}\fg_{n-p+i+1}(x))=0
\]
 which is equivalent, by Lemma 9, to
\[
\delta_{p+1}(x)+[m]_qq^{n-p+1}\sum^p_{i=0}\fg_{n-p+i+1}(x)=-\delta^{N-m+1}_{p+1}(y)\]
for some $y\in E^{n+m-N}$\\

On the other hand, one has by direct computation
\[
\delta_{p+1}(d^{N-m}_p y)+[m]_q
q^{n-p+1}\sum^p_{i=0}\fg_{n-p+i+1}(d^{N-m}_py)=\delta^{N-m+1}_{p+1}(y)
\]
which achieves the proof of the statement. $\Box$

\noindent {\bf Remark 3.} Let us define an endomorphism of degree one $d_{p,m}$
of $\oplus_{n\geq p-1}E^n$ by setting
$d_{p,m}=\delta_{p+1}+[m]_qq^{n-p+1}\sum^p_{i=0}\fg_{n-p+i+1}:E^n\rightarrow
E^{n+1}$. One has again $d^N_{p,m}=0$ and $d_{p,1}=d_p$ whereas
$d_{p,N-1}=d_{p+1}$.

\begin{corollary}
The homomorphisms $\Psi_\ast$ of Corollary 5 are isomomorphisms.
\end{corollary}

\noindent {\bf Proof.} Let $x\in E^{n+N}$ be such that $d^{N-1}_p x=0$. In view
of Lemma 10, this is equivalent to the existence of a $y\in E^{n+N-1}$ such
that $d_{p+1}(x-d_py)=0$ which implies the surjectivity of
$\Psi_\ast:H^{n+N}_{(1)}(E,d_{p+1})\rightarrow H^{n+N}_{(N-1)}(E,d_p)$. If
$x=d_pz$, then $d_{p+1}d_p(z-y)=\delta^2_{p+1}(z-y)=0$ implies
$z-y=\delta^{N-2}_{p+1}u$, in view of Lemma 9, and therefore
$d_p(z-y)=d_p\delta^{N-2}_{p+1}u=d^{N-1}_{p+1}u$ which implies that
$\Psi_\ast:H^{n+N}_{(1)}(E,d_{p+1})\rightarrow H^{n+N}_{(N-1)}(E,d_p)$ is also
injective and therefore is an isomorphism.\\

On the other hand let $x\in E^{n+N-1}$ be such that $d_px=0$. Then
$d_{p+1}d_px=\delta^2_{p+1}x=0$ implies $x=\delta^{N-2}_{p+1}x'$ with $x'\in
E^{n+1}$, in view of Lemma 9, and $d^{N-1}_{p+1}x'=0$ which implies the
surjectivity of $\Psi_\ast:H^{n+1}_{(N-1)}(E,d_{p+1})\rightarrow
H^{n+N-1}_{(1)}(E,d_p)$. If $x=d^{N-1}_py$ then $x=\delta^{N-2}_{p+1}d_{p+1}y$
which implies that $\Psi_\ast:H^{n+1}_{(N-1)}(E,d_{p+1})\rightarrow
H^{n+N-1}_{(1)}(E,d_p)$ is also injective and is therefore an isomorphism.\\

Finally the fact that $\Psi_\ast:H^0_{(1)}(E,d_1)\rightarrow
H^0_{(N-1)}(E,d_0)$ is an isomorphism follows from the injectivity of
$\delta^{N-2}_1:E^1\rightarrow E^{N-1}$. $\Box$\\

By using again Lemma 9 and Lemma 10 one deduces easily the following
improvement of Corollary~7.

\begin{corollary} The homomorphisms $\bar\Psi_{\ast p}$ $(p\in\mathbb N)$ and
the homomorphisms $\bar\Psi_{\ast p,\ell}$  $(p\in \mathbb N, \ell\in \{
1,\dots,N-2\})$  are isomorphisms.
\end{corollary}

\noindent{\bf Proof.} The only nontrivial isomorphisms which are not contained
in Corollary 7 are:\\
$\bar\Psi_{\ast p}:H^p(E)\stackrel{\simeq}{\rightarrow} H^p_{(N-1)}(E_p,d_p)$\\
$\bar\Psi_{\ast
p,\ell}:0\stackrel{\simeq}{\rightarrow}Z^{p-1+\ell}_{(1)}(E_p,d_p)$\\
$\bar\Psi_{\ast
p,\ell}:Z^{p+\ell}_{(1)}(E_{p+1},d_{p+1})\stackrel{\simeq}{\rightarrow}
H^{p+\ell}_{(N-1)}(E_p,d_p)$\\
for $\ell\in \{ 1,\dots,N-2\}$.\\

Let $x\in E^p$ be such that $d^{N-1}_px=0$. Then by Lemma 10, there is a $y\in
E^{p-1}$ for which $d_{p+1}(x-d_py)=d(x-dy)=0$ which implies that
$\bar\Psi_{\ast p}:H^p(E)\rightarrow H^p_{(N-1)}(E_p,d_p)$ is surjective and
therefore bijective since the injectivity is a tautology in view of $d=d_p$ on
$E^{p-1}$.\\
Let $x\in E^{p-1+\ell}$ with $1\leq\ell\leq N-2$ be such that $d_p x=0$. Then
$d_{p+1}d_px=\delta^2_{p+1}x=0$ which implies $x=0$ in view of Lemma 9 and
$p\leq p-1+\ell<p+N-2$. Finally let $x\in E^{p+\ell}$ with $1\leq\ell\leq N-2$
be such that $d^{N-1}_px=0$. Then, again by Lemma 10, there is a $y\in
E^{p+\ell-1}$ for which $d_{p+1}(x-d_py)=0$ which implies (by changing $p$ in
$p+1$) $x=d_py$ in view of above. This completes the proof of Corollary
8.~$\Box$\\

\noindent {\bf Proof of Theorem 2.} Corollary 8 implies that the homomorphisms
$\Phi_{\ast 1,p}$ of Corollary 6 are isomorphisms and that one has:
\[
H^{Nr+p-1+m}_{(1)}(E_p,d_p)=H^{Nr+p+m}_{(N-1)}(E_p,d_p)=0,\ \ \forall
m\in\{1,\dots,N-2\}.
\]
This implies, in view of Corollary 2 ($(ii)\Leftrightarrow (iii)$) that the
$\Delta_\ast$ induce isomorphisms implying that the $\Phi_{\ast m,p}$ are
isomorphisms $\forall m\in\{1,\dots,N-1\}.$~$\Box$\\

By using the equivalence $(iii)\Leftrightarrow (i)$ of Corollary 2, one can now
completely describe the homologies of the $N$-complexes $(E_p,d_p)$ in terms of
the (co)homology of the cosimplicial module $(E^n)$.

\begin{theorem}
One has the following isomorphisms\\
$H^{Nr+p-1}_{(m)}(E_p,d_p)\simeq H^{2r+p-1}(E)$ for $r\in \mathbb N$ with
$r\geq 1$\\
$H^{N(r+1)-m+p-1}_{(m)}(E_p,d_p)\simeq H^{2r+p}(E)$ for $r\in \mathbb N$ \\
$H^{p-1}_{(m)}(E_p,d_p)\simeq \ker (d:E^{p-1}\rightarrow E^p)$\\
and $H^n_{(m)}(E_p,d_p)=0$ otherwise (i.e. if $n\not= p-1 \mod (N)$
and\linebreak[4] $n+m\not= p-1 \mod (N))$.
\end{theorem}

Notice however that Theorem 2 is stronger than Theorem 3 in that it gives
explicitely the isomorphisms.\\

\noindent{\bf Remark 4.} The above results hold more generally under Assumption
$(A_0)$ for $\bk$ and $q\in \bk$ (i.e. $[N]_q=0$) for a pre-cosimplicial module
$E$ if we add the assumption $H_{(m)}(E_{p+1},\delta_{p+1})=0, \forall p\in
\mathbb N$ and $\forall m\in\{ 1,\dots,N-1\}$. However this latter assumption
is rather technical; it is why we have chosen to deduce it from Assumption
$(A_1)$ in the cosimplicial framework (Lemma 9).\\

In the case of a simplicial module ($E'_n$), one can proceed similarily by
introducing $N$-differentials $d'_p$ and $\delta'_p$ of degrees $-1$ and one
obtains similar results; the only difference being the replacement of the
indices $(m)$ by $(N-m)$ and of course $E'_{p-1}/d'(E'_p)$ instead of $\ker
(d:E^{p-1}\rightarrow E^p)$, (these replacements correspond to the duality).\\

The $N$-differential $d_0$ as well as its dual $d'_0$ have been considered in
several contexts \cite{May}, \cite{Kapr}, \cite{D-VK}, \cite{MD-V}, \cite{KW};
it is very natural. The $N$-differential $d_1$ was considered in \cite{D-VK},
\cite{MD-V}; in situations where one has products it is also very natural since
it satisfies a $q$-twisted version of the graded Leibniz rule (see below). In
these cases $(p=0,1)$, $E_p=E$ and one does not need to distinguish the degree
$p-1$. Let us summarize the result (i.e. Theorem 3 and its dual version) for
these cases.

\begin{theorem}
Let $\bk$ and $q\in\bk$ satisfy Assumption $(A_1)$, let $E$ be a cosimplicial
module and let $E'$ be a simplicial module. Then one has:
\[
\begin{array}{ll}
(0)& H^{Nr-1}_{(m)}(E,d_0)=H^{2r-1}(E),\ H^{N(r+1)-m-1}_{(m)}(E,d_0)=H^{2r}(E)\
\mbox{and}\\
& H^n_{(m)}(E,d_0)=0\ \ \mbox{otherwise},\\
\\
(0') & H_{(m),Nr-1}(E',d'_0)=H_{2r-1}(E'),H_{(m),Nr+m-1}(E',d'_0)=H_{2r}(E')\ \
\mbox{and}\\
& H_{(m),n}(E',d'_0)=0\ \ \mbox{otherwise},\\
\\
(1) & H^{Nr}_{(m)}(E,d_1)=H^{2r}(E), H^{N(r+1)-m}_{(m)}(E,d_1)=H^{2r+1}(E)\ \
\mbox{and}\\
& H^n_{(m)}(E,d_1)=0\ \ \mbox{otherwise},\\
\\
(1') & H_{(m),Nr}(E',d'_1)=H_{2r}(E'), H_{(m),Nr+m}(E',d'_1)=H_{2r+1}(E')\ \
\mbox{and}\\
& H_{(m),n}(E',d'_1)=0\ \ \mbox{otherwise},
\end{array}
\]
for $r\in\mathbb N$ and $m\in\{1,\dots,N-1\}$.

\end{theorem}

Let $N$ be a prime number, then $\bk=\mathbb Z/N\mathbb Z$ and $q=1\in \mathbb
Z/N\mathbb Z$ satisfy Assumption $(A_1)$. By applying the statement $(0')$ of
Theorem 4 to the simplicial module over $\mathbb Z/N\mathbb Z$ associated to a
simplicial complex (with $q=1$), one obtains an improvement of the results of
\cite{May}. In fact, applied to this situation, Statement $(0')$ of Theorem 4
implies the results of \cite{May} and the complete identification of the
generalized homology in terms of the ordinary homology with coefficients in
$\mathbb Z/N\mathbb Z$ of the simplicial complex.\\

Let $\bk$ and $q\in \bk$ satisfy Assumption $(A_1)$, let $\cala$ be an
associative unital $\bk$-algebra and let $\calm$ be an $\cala-\cala$ bimodule.
It is well known \cite{J-LL}, \cite{Weib} that the complex of $\calm$-valued
Hochschild cochains of $\cala$ is the complex of a cosimplicial module. Applied
to this case, Statement $(1)$ of Theorem 4 is exactly the result announced in
\cite{MD-V} expressing the generalized cohomology in terms of the Hochschild
cohomology. Dualy, the complex of Hochschild chains of $\cala$ with
coefficients in $\calm$ is the complex of a simplicial module and Theorem 1 of
\cite{KW} is Statement $(0')$ of Theorem 4 applied to this case with
$\calm=\cala^\ast$.\\

Borrowing a sentence from Kassel and Wambst, (in the introduction of
\cite{KW}), we can say that Theorem 3 and its dual version provide alternative
exotical definitions of the (co)homology of a (co)simplicial module, (whenever
there is a $q\in\bk$ such that Assumption $(A_1)$ is satisfied).

\newpage

\section{Generalization of graded differential algebras}

Throughout this section $q$ is a fixed element of $\bk$ such that Assumption
$(A_1)$ is satisfied, i.e. $[N]_q=0$ and $[n]_q$ is invertible in $\bk$ for any
$n\in\{1,\dots,N-1\}$. Let us introduce the following definition \cite{D-VK},
\cite{MD-V}. A {\sl graded q-differential algebra} is a $\mathbb N$-graded
associative unital $\bk$-algebra $\fraca=\oplus_{n\in\mathbb N}\fraca^n$
equipped with a $\bk$-linear endomorphism $d$ of degree 1, its $q$-{\sl
differential}, satisfying $d^N=0$ and the {\sl graded $q$-Leibniz rule}:
\[
d(\alpha\beta)=d(\alpha)\beta + q^a\alpha d(\beta),\ \ \forall \alpha\in
\fraca^a\ \ \mbox{and}\ \ \forall\beta\in \fraca.
\]
Thus, $(\fraca,d)$ is in particular a positive (cochain) $N$-complex.\\
This is clearly a generalization of the notion of graded differential algebra
since for $N=2$ and $q=-1$ it reduces to it, (for $N=2$, $\bk$ and $-1\in \bk$
obviously satisfy Assumption $(A_1)$).\\

The graded $q$-Leibniz rule implies that $Z_{(1)}(\fraca)$ is a graded unital
subalgebra of $\fraca$ and that $B_{(1)}(\fraca)$ is a graded two-sided ideal
of $Z_{(1)}(\fraca)$. Thus $H_{(1)}(\fraca)$ is a $\mathbb N$-graded
associative unital $\bk$-algebra.\\

One defines in an obvious manner the notions of homomorphisms of graded
$q$-differential algebras, of graded $q$-differential subalgebras, etc.\\

In the present context, graduations over $\mathbb Z_N=\mathbb Z/N\mathbb Z$ are
also very natural, (see e.g. Example 1 below), therefore let us define a {\sl
$\mathbb Z_N$-graded $q$-differential algebra} to be a $\mathbb Z_N$-graded
associative unital $\bk$-algebra equipped with a $q$-differential satisfying
the above axioms. A graded $q$-differential algebra as defined before can be
considered as a $\mathbb Z_N$-graded $q$-differential algebra by retaining only
the degree modulo $N$ and we speak then of the {\sl underlying $\mathbb
Z_N$-graded $q$-differential algebra} for the latter. On the other hand let
$\fraca=\oplus_{p\in \mathbb Z_N}\fraca^p$ be a $\mathbb Z_N$-graded algebra
and let $n\mapsto p(n)$ be the canonical projection of $\mathbb N$ onto
$\mathbb Z_N$. One can associate to $\fraca$ a $\mathbb N$-graded algebra
$p^\ast\fraca=\oplus_{n\in \mathbb N}\  p^\ast \fraca^n$ in the following
manner. An homogeneous element of $p^\ast\fraca$ is a pair $(n,\alpha)\in
\mathbb N \times \fraca^p$ with $p(n)=p\in \mathbb Z_N$ and we identify
$p^\ast\fraca^n=(n,\fraca^{p(n)})$ with the $\bk$-module $\fraca^{p(n)}$, the
product in $p^\ast\fraca$ being defined by
$(m,\alpha)(n,\beta)=(m+n,\alpha\beta)$. The canonical projection $\pi:p^\ast
\fraca\rightarrow \fraca$ defined by $\pi(n,\alpha)=\alpha$ is an homomorphism
of the $\mathbb Z_N$-graded algebra underlying to $p^\ast\fraca$ onto $\fraca$
and the pair $(p^\ast\fraca,\pi)$ is characterized by the following universal
property: For any $\mathbb N$-graded algebra $\Omega$ and for any homomorphism
$\varphi$ of $\mathbb Z_N$-graded algebras of the underlying $\mathbb
Z_N$-graded algebra to $\Omega$ into $\fraca$ there is a unique homomorphism of
$\mathbb N$-graded algebras $\bar\varphi:\Omega\rightarrow p^\ast\fraca$ such
that $\varphi=\pi \circ \bar \varphi$. If $\fraca$ is a $\mathbb Z_N$-graded
$q$-differential algebra then there is a unique $q$-differential on
$p^\ast\fraca$ such that $\pi$ is an homomorphims of the underlying $\mathbb
Z_N$-graded $q$-differential algebra of $p^\ast\fraca$ onto $\fraca$.\\

In what follows we give some examples of constructions of graded
$q$-differential algebras.
\newpage

\noindent {\bf Construction 1: $q$-derivations.}\\

Let $\Bg=\oplus_{n\in \mathbb N}\Bg^n$ be a $\mathbb N$-graded associative
unital $\bk$-algebra and let $d$ be a $\bk$-linear mapping of degree 1 of $\Bg$
into itself which satisfies the graded $q$-Leibniz rule; we shall refer to such
a $d$ as a $q$-{\sl derivation of degree} 1 of $\Bg$. Then, by using the
inductive definition of the $q$-binomial coefficients $\left[\begin{array}{l}
n\\m
\end{array}\right]_q$ (see in Section 1), one obtains by induction on $n$:
\[
d^n(\alpha\beta)=\sum^n_{m=0}q^{a(n-m)}\left[\begin{array}{l} n\\m
\end{array}\right]_qd^m(\alpha)d^{n-m}(\beta),\ \forall \alpha\in \Bg^a\
\mbox{and}\ \forall \beta\in\Bg.
\]
On the other hand it follows from $[N]_q=0$ and from the invertibility of
$[m]_q$ for $m\in\{1,\dots,N-1\}$ that one has $\left[\begin{array}{l} N\\m
\end{array}\right]_q=0$,  $\forall m\in \{1,\dots,N-1\}$. This implies that
$d^N$ is a homogeneous derivation of $\Bg$ (i.e.
$d^N(\alpha\beta)=d^N(\alpha)\beta +\alpha d^N(\beta))$. It follows that
$\ker(d^N)$ is a graded unital subalgebra of $\Bg$ which is stable by $d$ and
therefore that $\fraca=\ker(d^N)$ equipped with $d$ is a graded
$q$-differential algebra. This is in fact the very reason of the compatibility
of $d^N=0$ with the graded $q$-Leibniz rule for $d$ under Assumption $(A_1)$.
This provides the construction of a graded $q$-differential algebra by starting
from a $q$-derivation of degree 1.\\

There is an easy way to produce such a $q$-derivation of degree 1 of $\Bg$ by
starting from an element $e$ of $\Bg^1$: Define the $\bk$-linear mapping
$ad_q(e)$ of $\Bg$ into itself by setting $ad_q(e)(\alpha)=e\alpha-q^a\alpha
e$,  $\forall \alpha\in \Bg^a$; then it is immediate that $ad_q(e)$ satisfies
the graded $q$-Leibniz rule. The $q$-derivations of degree 1 obtained by this
construction will be refered to as the {\sl inner $q$-derivations} of degree 1
of $\Bg$. By induction on $n$, one obtains:
\[
(ad_q(e))^n(\alpha)=\sum^n_{m=0}(-1)^m
q^{ma+\frac{m(m-1)}{2}}\left[\begin{array}{l}n\\ m\end{array}\right]_q
e^{n-m}\alpha e^m,\ \ \forall \alpha\in \Bg^a
\]
In view of Assumption $(A_1)$, $\left[\begin{array}{l}N\\
m\end{array}\right]_q=0$ for $m\in\{1,\dots,N-1\}$ so one has
$ad_q(e)^N(\alpha)=e^N\alpha+ (-1)^Nq^{\frac{N(N-1)}{2}}\alpha e^N$. By using
again Assumption $(A_1)$ one can show that $(-1)^Nq^{\frac{N(N-1)}{2}}=-1$;
this is obvious if $N$ is odd and for $N$ even, this follows from the
invertibility of $[N/2]_q$ which implies $q^{N/2}=-1$. Thus one has
$(ad_q(e))^N=ad(e^N)$, i.e. $(ad_q(e))^N(\alpha)=e^N\alpha-\alpha e^N$ and
therefore $\ker((ad_q(e))^N)$ is the commutant $\{e^N\}'$ of $e^N$ in $\Bg$. In
particular, if $e$ is an element of degree 1 of $\fraca$ such that $e^N$ is
central, then $ad_q(e)$ is a $q$-differential on $\fraca$, i.e. $\fraca$
equipped with $d=ad_q(e)$ is a graded $q$-differential algebra.
Before giving an example of this construction, let us generalize the notion
$q$-derivation of degree 1.\\

A {\sl $q$-derivation of degree} $\ell \in \mathbb Z$ of $\Bg$ is an
homogeneous $\bk$-linear mapping\linebreak[4] $L:\Bg \rightarrow \Bg$ of degree
$\ell$  such that :
\[
L(\alpha\beta)=L(\alpha)\beta + q^{\ell a}\alpha L(\beta),\ \ \forall \alpha\in
\Bg^a\ \mbox{and}\ \forall \beta\in \Bg.
\]
The module of all $q$-derivations of degree $\ell$ of $\Bg$ will be denoted by
$\gder^\ell_q(\Bg)$ and the direct sum $\gder_q(\Bg)=\oplus_{\ell\in \mathbb Z}
\gder^\ell_q(\Bg)$ will be refered to as the module of {\sl graded
$q$-derivations} of $\Bg$. Again, there is an easy way to produce a
$q$-derivation of degree $\ell\in \mathbb N$ by starting from an element
$\lambda\in \Bg^\ell$: Define the $\bk$-linear mapping $ad_q(\lambda)$ of $\Bg$
into itself by setting $ad_q(\lambda)(\alpha)=\lambda \alpha-q^{\ell a}\alpha
\lambda,\ \ \forall \alpha\in \Bg^a$; one easily verifies that $ad_q(\lambda)$
is a $q$-derivation of degree $\ell$ of $\Bg$. The direct sum
$Int_q(\Bg)=\oplus_\ell ad_q(\Bg^\ell)$ is a graded submodule of $\gder_q(\Bg)$
which will be refered to as the module of {\sl inner graded $q$-derivations}.
Notice that there is no natural graded bilinear bracket on $\gder_q(\Bg)$ for
$q\not=1$ and $q\not=-1$. This does not mean that there are no multilinear
operation, for instance the above formulae imply that the symmetrized
composition of $N$ $q$-derivations of degree one is a $q$-derivation of degree
$N$.\\
All this obviously also applies to the $\mathbb Z_N$-graded case.\\

\noindent{\bf Example 1: Matrix algebra.}\\

Let us introduce the standard basis $E^k_\ell,\ \ (k,\ell\in \{1,\dots,N\})$,
of the algebra $M_N(\bk)$ of $N\times N$ matrices defined by
$(E^k_\ell)^i_j=\delta^k_j\delta^i_\ell$. One has $E^k_\ell
E^r_s=\delta^k_sE^r_\ell$ and $\sum^N_{n=1}E^n_n=\bbbone$. It follows that one
can equip $M_N(\bk)$ with a structure of\linebreak[4] $\mathbb Z_N$-graded
algebra by giving to $E^k_\ell$ the degree $k-\ell$ mod($N$).\\
Let $e=\lambda_1E^2_1+\dots +\lambda_{N-1}E^N_{N-1}+\lambda_N E^1_N$,
($\lambda_i\in \bk$), be an arbitrary element of degree one of $M_N(\bk)$. One
has $e^N=\lambda_1\dots \lambda_N\bbbone$ and therefore $e^N$ is in the center
of $M_N(\bk)$ so one defines a $q$-differential $d$ on $M_N(k)$ by setting
$d(A)=eA-q^a Ae$ for $A\in M_N(\bk)^a$. Equipped with the above structures
$M_N(\bk)$ is a $\mathbb Z_N$-graded $q$-differential algebra and therefore
$p^\ast M_N(\bk)$ is a graded $q$-differential algebra.\\

\noindent{\bf Construction 2: (pre-)cosimplicial algebras.}\\

Let $(E^n)_{n\in \mathbb N}$ be a pre-cosimplicial module (over $\bk$) with
coface homomorphisms $\fg_i:E^n\rightarrow E^{n+1}$ as before. We shall say
that $(E^n)$ is a {\sl pre-cosimplicial algebra} if the $\mathbb N$-graded
module $E=\oplus_{n\in\mathbb N} E^n$ is equipped with a structure of $\mathbb
N$-graded associative unital $\bk$-algebra compatible with its structure of
graded module and such that one has $\forall \alpha\in E^a$, $\forall \beta\in
E^b$\\

\noindent $(\fraca\Fg_1)\hspace{1.5cm} \fg_i(\alpha\beta)=\left\{
\begin{array}{lll}
\fg_i(\alpha)\beta &\mbox{if} & i\leq a\\
\alpha \fg_{i-a}(\beta) &\mbox{if} & i>a
\end{array}
\right.
, i\in\{0,\dots,a+b+1\}$\\
and\\
$(\fraca\Fg_2)\hspace{1.5cm} \fg_{a+1}(\alpha)\beta = \alpha \fg_0(\beta)$\\
where $(\alpha,\beta)\mapsto \alpha\beta$ denote the product of $E$. If
furthermore $(E^n)$ is a cosimplicial module (i.e. if one has the codegeneracy
homomorphisms) then $(E^n)$ will be called a {\sl cosimplicial algebra}
whenever one has $\forall \alpha\in E^a$, $\forall \beta\in E^b$\\

\noindent $(\fraca\Sg)\hspace{1.5cm} \sg_i(\alpha\beta)=\left\{
\begin{array}{lll}
\sg_i(\alpha)\beta &\mbox{if} & i< a\\
\alpha \sg_{i-a}(\beta) &\mbox{if} & i\geq a
\end{array}
\right.
, i\in\{0,\dots,a+b-1\}$\\

Let the $N$-differential $d_1$ of $E$ be defined as in Section 3,
i.e.\linebreak[4] $d_1=\sum^n_{i=0}q^i\fg_i-q^n\fg_{n+1}:E^n\rightarrow
E^{n+1}$. Then one has the following result.

\begin{proposition}
Let $(E^n)$ be a pre-cosimplicial algebra. Then $E=\oplus_nE^n$ equipped with
$d_1$ is a graded $q$-differential algebra.
\end{proposition}

\noindent{\bf Proof.} Since, from Corollary 4, we know that $d^N_1=0$ it is
sufficient to prove that $d_1$ satisfies the graded $q$-Leibniz rule. Let
$\alpha\in E^a$ and $\beta\in E^b$, one has
$d_1(\alpha\beta)=\sum^{a+b}_{i=0}q^i\fg_i(\alpha\beta)-q^{a+b}\fg_{a+b+1}(\alpha\beta)$ which can be rewritten by applying $(\fraca\Fg_1)$ as\\
$d_1(\alpha\beta)=\sum^a_{i=0}q^i\fg_i(\alpha)\beta +
q^a\alpha\left(\sum^b_{j=1}q^j\fg_j(\beta)-q^b\fg_{b+1}(\beta)\right)$
and, by using $(\fraca\Fg_2)$\\
$d_1(\alpha\beta)=\left(\sum^a_{i=0}q^i\fg_i(\alpha)-q^a\fg_{a+1}(\alpha)\right)\beta+q^a\alpha\left(\sum^b_{j=0}q^j\fg_j(\beta)-q^b\fg_{b+1}(\beta)\right)$ i.e.,\\
$d_1(\alpha\beta)=d_1(\alpha)\beta + q^a\alpha d_1(\beta).\Box$\\

Applied to the case $N=2$, $q=-1$, Proposition 3 implies that $E$ equipped with
the usual $d$ is a graded differential algebra, which implies in particular
that the usual cohomology $H(E)$ of $(E^n)$ is a $\mathbb N$-graded associative
unital $\bk$-algebra.\\

In the case where $(E^n)$ is a cosimplicial algebra, Theorem 4 (1) implies that
$H^{Nr}_{(m)}(E,d_1)=H^{2r}(E)$, $H^{N(r+1)-m}_{(m)}(E,d_1)=H^{2r+1}(E)$ and
$H^n_{(m)}(E,d_1)=0$ otherwise.\\

\noindent{\bf Warning.} $H_{(1)}(E,d_1)$ is an algebra (as shown above) as well
as $H(E)$. However, in spite of the fact that the module $H_{(1)}(E,d_1)$ does
only depend on the module $H(E)$ in the cosimplicial case, the algebra
structures are not the same e.g.\linebreak[4]
$H^{N(r+1)-1}_{(1)}(E,d_1)=H^{2r+1}(E)$ and
$H^{N(s+1)-1}_{(1)}(E,d_1)=H^{2s+1}(E)$ but\linebreak[4]
$H^{N(r+1)-1}_{(1)}H^{N(s+1)-1}_{(1)}\subset H^{N(r+s+2)-2}_{(1)}=0$ if $N\geq
3$ although $H^{2r+1}(E)H^{2s+1}(E)$ can be $\not= 0$.\\

\noindent{\bf Example 2: Simplicial forms.}\\

Let $K$ be a simplicial complex, i.e. $K$ is a set equipped with a set $\cals$
of non-empty subsets such that $X\in \cals$ and $Y\subset X$ with
$Y\not=\emptyset$ implies $Y\in \cals$ and such that $\{x\}\in \cals$, $\forall
x\in K$. An {\sl ordered n-simplex} ($n\in \mathbb N$) is a sequence
$(x_0,\dots,x_n)$ with $x_i\in K$ such that $\im(x)=\cup^n_{i=0}\{x_i\}\in
\cals$. A $\bk$-{\sl valued  simplicial $n$-form} is a $\bk$-valued function
$(x_0,\dots,x_n)\mapsto \omega(x_0,\dots,x_n)$ on the set of ordered
$n$-simplices. Let $\Omega^n(K,\bk)$ denote the $\bk$-module of all
$\bk$-valued simplicial $n$-forms. The graded module
$\Omega(K,\bk)=\oplus_{n\in \mathbb N}\Omega^n(K,\bk)$ is a $\mathbb N$-graded
associative $\bk$-algebra for the product $(\alpha,\beta)\mapsto \alpha\beta$
defined by
$(\alpha\beta)(x_0,\dots,x_{a+b})=\alpha(x_0,\dots,x_a)\beta(x_a,\dots,x_{a+b})$ for $\alpha\in\Omega^a(K,\bk)$, $\beta\in\Omega^b(K,\bk)$ and any ordered $(a+b)$-simplex $(x_0,\dots,x_{a+b})$. One verifies easily that $(\Omega^n(K,\bk))$ is a cosimplicial algebra with $\fg_i$, $\sg_i$ defined by $\fg_i(\omega)(x_0,\dots,x_{n+1})=\omega(x_0,\oisi,x_{n+1})$ $\forall i\in\{0,\dots,n+1\}$ (where $\omis$   means omission of $x_i$) and $\sg_i(\omega)(x_0,\dots,x_{n-1})=\omega(x_0,\dots,x_i,x_i,\dots,x_{n-1})$ $\forall i\in \{0,\dots,n-1\}$ for $\omega\in \Omega^n(K,\bk)$ and any ordered $(n+1)$-simplex $(x_0,\dots,x_{n+1})$ and any ordered $(n-1)$-simplex $(x_0,\dots,x_{n-1})$. It follows that equipped with $d_1$, $\Omega(K,\bk)$ is a graded $q$-differential algebra. Notice that $d_1$ is then given by :
\[
d_1(\omega)(x_0,\dots,x_{n+1})=\sum^n_{i=o}q^i\omega(x_0,\oisi,x_{n+1})-q^n\omega(x_0,\dots,x_n)
\]
for $\omega\in \Omega^n(K,\bk)$ and any ordered $(n+1)$-simplex
$(x_0,\dots,x_{n+1})$.\\
According to Theorem 4 (1) one has :
\[
H^{Nr}_{(m)}(\Omega(K,\bk),d_1)=H^{2r}(K,\bk),\ \  H^{N(r+1)-m}_{(m)}(\Omega
(K,\bk),d_1)=H^{2r+1}(K,\bk)
\]
and $H^n_{(m)}(\Omega(K,\bk),d_1)=0$ otherwise, where $H(K,\bk)$ denotes the
usual cohomology of the simplicial complex with coefficients in $\bk$.\\

\noindent{\bf Example 3: Hochschild cochains.}\\

Let $\cala$ be an associative unital $\bk$-algebra and let $\calm$ be an
$\cala-\cala$ bimodule. Recall that a $\calm$-{\sl valued  Hochschild cochain}
of degree $n$ ($n\in \mathbb N$) is a $n$-linear mapping of
$\underbrace{\cala\times\dots\times\cala}_n$ into $\calm$. The $\bk$-module of
all $\calm$-valued Hochschild cochains is denoted by $C^n(\cala,\calm)$. The
sequence $(C^n(\cala,\calm))$ is a cosimplicial module with cofaces $\fg_i$ and
codegeneracies $\sg_i$ defined by \cite{J-LL}, \cite{Weib}\\
$\fg_0(\omega)(x_0,\dots,x_n)=x_0\omega(x_1,\dots,x_n)$\\
$\fg_i(\omega)(x_0,\dots,x_n)=\omega(x_0,\dots,x_{i-1}x_i,\dots,x_n)$\hspace{2cm} for $i\in\{1,\dots,n\}$\\
$\fg_{n+1}(\omega)(x_0,\dots,x_n)=\omega(x_0,\dots,x_{n-1})x_n$\\
and\\
$\sg_i(\omega)(x_1,\dots,x_{n-1})=\omega(x_1,\dots,x_i,\bbbone,x_{i+1},\dots,x_{n-1})$ \hspace{1cm} for $i\in \{0,\dots,n-1\}$\\
for $\omega\in C^n(\cala,\calm)$ and $x_i\in \cala$.\\

Thus $C(\cala,\calm)=\oplus_n C^n(\cala,\calm)$ equipped with $d_1$ is a
$N$-complex and one has $H^{Nr}_{(m)}(C(\cala,\calm),d_1)=H^{2r}(\cala,\calm),\
\ H^{N(r+1)-m}_{(m)}(C(\cala,\calm),d_1)=H^{2r+1}(\cala,\calm)$ and
$H^n_{(m)}(C(\cala,\calm),d_1)=0$ otherwise, where $H(\cala,\calm)$ denote the
usual Hochschild cohomology of $\cala$ with values in $\calm$. The
$N$-differential $d_1$ being given by
\[
\begin{array}{ll}
d_1(\omega)(x_0,\dots,x_n)& =x_0\omega(x_1,\dots,x_n)
\\
\\
&+\sum^n_{i=1}q^i\omega(x_0,\dots,x_{i-1}x_i,\dots,x_n)-q^n\omega(x_0,\dots,x_{n-1})x_n
\end{array}
\]
for $\omega\in C^n(\cala,\calm)$ and $x_i\in \cala$. This $N$-differential has
been introduced in \cite{D-VK} with another notation, $(d_1\mapsto
\delta_q)$.\\

In the case $\calm=\cala$, $C(\cala,\cala)$ has a natural structure of $\mathbb
N$-graded associative unital $\bk$-algebra with product $(\alpha,\beta)\mapsto
\alpha\beta$ given by
\[
\alpha\beta(x_1,\dots,x_{a+b})=\alpha(x_1,\dots,x_a)\beta(x_{a+1},\dots,x_{a+b}),
\]
for $\alpha\in C^a(\cala,\cala),\ \beta\in C^b(\cala,\cala), x_i\in \cala$.\\
It is easily verified that $(C^n(\cala,\cala))$ is a cosimplicial algebra and
therefore it follows that $(C(\cala,\cala),d_1)$ is a graded $q$-differential
algebra.\\

\noindent{\bf Example 4: The tensor algebra over $\cala$ of
$\cala\otimes\cala$.}\\

Let again $\cala$ be an associative unital $\bk$-algebra and let us denote by
$\ft(\cala)=\oplus_{n\in \mathbb N}\ft^n(\cala)$ the tensor algebra over
$\cala$ of the $\cala-\cala$ bimodule $\cala\otimes \cala$, ($\otimes$ always
means the tensor product over $\bk$). This is a $\mathbb N$-graded associative
unital $\bk$-algebra with $\ft^n(\cala)=\otimes^{n+1}\cala$ and product
$(x_0\otimes \dots\otimes x_n)(y_0\otimes \dots \otimes
y_m)=x_0\otimes\dots\otimes x_{n-1}\otimes x_ny_0\otimes y_1\otimes \dots
\otimes y_m$ for $x_i,y_j\in \cala$. One defines a structure of cosimplicial
module for $(\ft^n(\cala))$ by setting\\
$\fg_0(x_0\otimes\dots \otimes x_n)=\bbbone \otimes x_0\otimes\dots\otimes
x_n$\\
$\fg_i(x_0\otimes\dots\otimes x_n)=x_0\otimes\dots \otimes x_{i-1}\otimes
\bbbone \otimes x_i\otimes\dots\otimes x_n\  \mbox{for}\ 1\leq i\leq n$\\
$\fg_{n+1}(x_o\otimes \dots\otimes x_n)=x_0\otimes\dots\otimes x_n\otimes
\bbbone$\\
and\\
$\sg_i(x_0\otimes\dots\otimes x_n)=x_0\otimes\dots\otimes x_ix_{i+1}\otimes
\dots\otimes x_n$ for $0\leq i\leq n-1$\\
One shows easily that with the above structures,  $\ft(\cala)$ is a
cosimplicial algebra and therefore, $\ft(\cala)$ equipped with the
$N$-differential $d_1$ is a graded $q$-differential algebra. It follows
(Theorem 4 (1)) that the $H^n_{(m)}(\ft(\cala),d_1)$ can be computed in terms
of the cohomology $H^n(\ft(\cala))$ of the cosimplicial module $\ft(\cala)$.
However we shall give later a direct proof of the triviality of these
generalized cohomologies whenever $\cala$ admits a linear form $\omega$
satisfying $\omega(\bbbone)=1$ by using Lemma 5.\\
In this case $d_1$ is given by
\[
d_1(x_0\otimes\dots\otimes x_n)=\sum^n_{i=0} q^i x_0\otimes\dots \otimes
x_{i-1}\otimes \bbbone\otimes x_i\otimes \dots\otimes x_n - q^n x_0\otimes\dots
\otimes x_n\otimes \bbbone
\]
so it coincides with the $q$-differential of $\ft(\cala)$ introduced in
\cite{D-VK}, \cite{MD-V}.\\

As pointed out in \cite{D-VK}, the graded algebra $\ft(\cala)$ can be
characterized by a universal property. Indeed $\cala\otimes \cala$ is the free
$\cala-\cala$ bimodule generated by $\tau=\bbbone\otimes \bbbone$, hence
$\ft(\cala)$ is the $\mathbb N$-graded associative $\bk$-algebra generated by
$\cala$ in degree 0 and by a free generator $\tau$ of degree 1;
$x_0\otimes\dots\otimes x_n=x_0\tau x_1\dots \tau x_n$. This implies the
following result \cite{D-VK}.

\begin{proposition}
Let $\fraca=\oplus_n\fraca^n$ be a $\mathbb N$-graded associative unital
$\bk$-algebra. Then for any homomorphism of unital $\bk$-algebras
$\varphi:\cala\rightarrow \fraca^0$ and for any element $\alpha$ of $\fraca^1$,
there is a unique homomorphism of graded unital $\bk$-algebras\linebreak[4]
$\ft_{\varphi,\alpha}:\ft(\cala)\rightarrow \fraca$ which extends $\varphi$ and
is such that
$\ft_{\varphi,\alpha}(\tau)=\alpha,\linebreak[4](\tau=\bbbone\otimes \bbbone
\in \ft^1(\cala))$.
\end{proposition}

By applying this result to the case where $\fraca=C(\cala,\cala)$, (Example 3),
where $\varphi$ is the identity mapping of $\cala$ onto itself considered as an
homomorphism of $\cala$ into $C^0(\cala,\cala)(=\cala)$ and where $\alpha$ is
again the identity mapping of $\cala$ onto itself considered as an element of
$C^1(\cala,\cala)$, one obtains the canonical homomorphism
$\Psi=\ft_{\varphi,\alpha}:\ft(\cala)\rightarrow C(\cala,\cala)$ of graded
unital algebras defined in \cite{Mass} which reads
\[
\Psi(x_0\otimes\dots\otimes x_n)(y_1,\dots,y_n)=x_0y_1x_1\dots y_nx_n, \ \ \
\mbox{for} x_i,y_j\in \cala.
\]

\begin{proposition}
The above homomorphism $\Psi$ is an homomorphism of cosimplicial algebras, i.e.
one has $\Psi\circ \fg_i=\fg_i\circ \Psi$ and $\Psi\circ \sg_i=\sg_i\circ\Psi$
with obvious notations.
\end{proposition}

This statement is easy to check. This implies in particular that one has
$\Psi\circ d_1=d_1\circ \Psi$ as well as $\Psi\circ d=d\circ \Psi$ where the
usual cosimplicial differential $d$ is the ordinary Hochschild differential on
$C(\cala,\cala)$.\\

As pointed out in \cite{D-VK} (with a slightly different formulation), the
$q$-differential $d_1$ on $\ft(\cala)$ is the unique $q$-derivation of degree 1
of $\ft(\cala)$ such that one has $d_1(x)=\bbbone \otimes x - x\otimes \bbbone$
for $x\in \cala$ and $d_1(\bbbone \otimes \bbbone)=\bbbone \otimes
\bbbone\otimes \bbbone (=(\bbbone \otimes \bbbone)^2)$.\\
By induction on the integer $n$, one shows that
\[
d^n_1(\bbbone\otimes\bbbone)=[n]_q!(\bbbone\otimes\bbbone)^{n+1}=[n]_q!\
\bbbone^{\otimes(n+2)}
\]
and
\[
d^n_1(x)=[n]_q!(\bbbone\otimes \bbbone)^{n-1}d(x)=[n]_q!\ \bbbone^{\otimes
n}d(x)
\]
for $x\in \cala$ with $d(x)=\bbbone\otimes x-x\otimes \bbbone=d_1(x)$.\\

\noindent {\bf Remark 5.} One can extend the setting of pre-cosimplicial
algebras in a framework of monoidal categories of pre-cosimplicial modules.
More precisely let $E$, $F$ and $G$ be 3 pre-cosimplicial modules with coface
homomorphisms denoted by $\fg_i$ and assume that one has bilinear mappings
$\cup:E^a\times F^b\rightarrow G^{a+b}$ such that

\[
\fg_i(\alpha\cup \beta)=\left\{\begin{array}{ll}
\fg_i(\alpha)\cup\beta & \mbox{if}\ \ i\leq a\\
\alpha\cup \fg_{i-a}(\beta) & \mbox{if}\ \ i>a
\end{array}\right., \ \ i\in \{0,\dots,a+b+1\}
\]
and such that\\
\[
\fg_{a+1}(\alpha)\cup\beta = \alpha\cup \fg_0(\beta),\ \ \forall \alpha\in E^a\
\ \mbox{and}\ \ \forall \beta\in F^b.
\]
Then one has $d_1(\alpha\cup\beta)=d_1(\alpha)\cup \beta+q^a\alpha \cup
d_1(\beta)$, with $d_1$ defined as in Section~3. This applies in particular to
a generalization of Example 2 by taking for E the simplicial forms with
coefficients in a $\bk$-module $\cal E$, for $F$ the simplicial forms with
coefficients in a $\bk$-module $\calf$ and for $G$ the simplicial forms with
coefficients in the $\bk$-module $\cal E \otimes \cal F$, $\cup$ being then the
tensor product over $\bk$ combined with the product of simplicial forms. This
also applies to Hochschild cochains by taking $E=C(\cala,\calm)$,
$F=C(\cala,\caln)$ and $G=C(\cala,\calm\otimes_{\cala}\caln)$ where $\calm$ and
$\caln$ are $\cala-\cala$ bimodules; the bilinear mapping $\cup$ is then the
tensor product over $\cala$ and corresponds to the usual cup product, (see in
\cite{D-VK}).\\

\noindent {\bf Construction 3: Universal $q$-differential envelopes.}\\

Let $\cala$ be an associative unital $\bk$-algebra and consider the following
category $q_{\cala}$. An object of $q_{\cala}$ is a graded $q$-differential
algebra $\Omega=\oplus_{n\in \mathbb N}\Omega^n$ together with an homomorphism
of unital algebras $\varphi:\cala \rightarrow \Omega^0$; a morphism of
$(\Omega, \varphi)$ into $(\Omega',\varphi')$ is an homomorphism
$\psi:\Omega\rightarrow \Omega'$ of graded $q$-differential algebras such that
$\varphi'=\psi\circ \varphi:\cala\rightarrow \Omega^{\prime 0}$. An object of
$q_{\cala}$ will be refered to as a $q$-{\sl differential calculus for}
$\cala$. An important property of $q_{\cala}$ is that it possesses an initial
universal object. Namely there exist a graded $q$-differential algebra
$\Omega_q(\cala)$ and an homomorphism of unital algebras $\varphi_q$ of $\cala$
into $\Omega^0_q(\cala)$ such that, for any $q$-differential calculus
$(\Omega,\varphi)$ for $\cala$, there is a unique homomorphism
$\psi_\Omega:\Omega_q(\cala)\rightarrow \Omega$ of graded $q$-differential
algebras satisfying $\varphi=\psi_\Omega\circ \varphi_q$. Clearly
$(\Omega_q(\cala),\varphi_q)$ is unique up to an isomorphism and will be
refered to as {\sl the universal $q$-differential calculus for} $\cala$ and the
graded $q$-differential algebra $\Omega_q(\cala)$ will be called {\sl the
universal $q$-differential envelope of $\cala$}. It is straightforward that
$\varphi_q$ is an isomorphism and we shall identify $\cala$ with
$\Omega^0_q(\cala)$. One constructs easily $\Omega_q(\cala)$ by generators and
relations \cite{D-VK}: it is the graded algebra generated by $\cala$ in degree
0 and by the $d^\ell(\cala)$ in degrees $\ell$ for $\ell\in \{1,\dots,N-1\}$
together with the relations of $\cala$ and the relations
\[
d^n(xy)=\sum^n_{m=0}\left[\begin{array}{c}
n\\ m
\end{array}
\right]_q d^m(x)d^{n-m}(y)\hspace{1cm} \forall x,y\in \cala
\]
and one defines the $d$ on $\Omega_q(\cala)$ by the graded $q$-Leibniz rule and
$d^N=0$ (from $d$ on $\cala$). One has $d(\bbbone)=0$ and the $\bk$-modules
$d^\ell(\cala)$ are all isomorphic to $\cala/\bk\bbbone$.\\

The graded $q$-differential algebra ($\ft(\cala),d_1$) of Example 4 together
with the identity $\cala=\ft^0(\cala)$ is a graded $q$-differential calculus
for $\cala$ and therefore there is a unique homomorphism of graded
$q$-differential algebras of $\Omega_q(\cala)$ into $\ft(\cala)$ which extends
the identity mapping of $\cala$ onto itself and the following result is not
very hard to prove \cite{D-VK}.

\begin{theorem}
The unique homomorphism of graded $q$-differential algebras of
$\Omega_q(\cala)$ into $\ft(\cala)$ equipped with $d_1$ which induces the
identity mapping of $\cala$ onto itself is injective.
\end{theorem}

It follows that we can identify $\Omega_q(\cala)$ with the graded
$q$-differential subalgebra of ($\ft(\cala),d_1$) generated by $\cala(\subset
\ft(\cala))$.

\begin{lemma}
Let $\cala^\ast$ be the dual module of $\cala$. The following conditions $(i)$
and $(ii)$ are equivalent for the algebra $\cala$.\\
$(i)$ The canonical bilinear mapping of $\cala^\ast\times \cala$ into $\bk$ is
surjective.\\
$(ii)$ There exists $\omega\in \cala^\ast$ such that $\omega(\bbbone)=1$.
\end{lemma}

\noindent {\bf Proof.} $(ii) \Rightarrow (i)$ is obvious since the canonical
image of $\cala^\ast\times \cala$ is an ideal of $\bk$. Assume that $(i)$ is
satisfied so that there exist $\omega_0\in \cala^\ast$ and $x_0\in \cala$ such
that $\omega_0(x_0)=1$; then $\omega\in \cala^\ast$ defined by
$\omega(x)=\omega_0(xx_0)$\ \  $\forall x\in \cala$ satisfies $(i)$. $\Box$

\begin{proposition}
Assume that $\cala$ satisfies the equivalent condition of Lemma~11. Then the
generalized (co)homologies of $(\ft(\cala),d_1)$ and of $\Omega_q(\cala)$ are
trivial in the sense that one has :
\[
\begin{array}{ll}
H^n_{(m)}(\ft(\cala),d_1)=H^n_{(m)}(\Omega_q(\cala))=0 & for\ \ n\geq 1\ \
and\\
\\
H^0_{(m)}(\ft(\cala),d_1)=H^0_{(m)}(\Omega_q(\cala))=\bk & for\ \ m\in
\{1,\dots,N-1\}.
\end{array}
\]
\end{proposition}

\noindent{\bf Proof.} Consider the $\mathbb Z$-graded module $E=\bk
e_{-(N-1)}\oplus \dots \oplus \bk e_{-1}\oplus\ft(\cala)$, where $\bk e_i$ is
the free $\bk$-module of rank 1 with generator $e_i$, $i\in
\{-(N-1),\dots,-1\}$. One extends $d_1$ to $E$ by setting $\tilde d_1
e_{-1}=\bbbone$ and $\tilde d_1 e_{-i}=e_{-(i-1)}$ for $N-1\geq i \geq 2$. One
still has $\tilde d_1^N=0$ so $(E,\tilde d_1)$ is a $N$-complex. Let  $\omega$
be a linear form on $\cala$ satisfying $\omega(\bbbone)=1$ and let us define an
endomorphism $h$ of degree $-1$ of $E$ by:
\[
h(x_0\otimes \dots \otimes x_n)=\omega(x_0)x_1\otimes \dots \otimes x_n
\]
for $x_i\in \cala$ and $n\geq 1$, $h(x_0)=-q^{-1}\omega(x_0)e_{-1}$ for $x_0\in
\cala$.
\[
h(e_{-i})=-q^{-(i+1)}(1+q+\dots+q^i)e_{-(i+1)}
\]
for $N-2\geq i\geq 1$  and $h(e_{-(N-1)})=0$. One has on $E$:  $h\tilde
d_1-q\tilde d_1 h=I$. By Lemma 5 it follows that $H_{(m)}(E)=0$, $\forall m\in
\{1,\dots,N-1\}$. The submodule $F=\bk e_{-(N-1)}\oplus \dots \oplus \bk
e_{-1}\oplus \Omega_q(\cala)$ is stable by $\tilde d_1$ and by $h$ and
therefore $H_{(m)}(F)=0$, $\forall m\in \{ 1,\dots,N-1\}$. This implies
immediately the result. The stability of $F$ under the action of $h$ follows
from the fact that $\Omega_q(\cala)$ is generated, as $\bk$-module, by the
$x_0d^{n_1}x_1\dots d^{n_p} x_p$, $x_i\in \cala$, $n_i\in \{1,\dots,N-1\}$ and
from the expression of the $d^nx$ as $d^nx=[n]_q!\bbbone^{\otimes n} dx$ for
$x\in\cala$. $\Box$\\

\noindent It is worth noticing here that the equivalent conditions of Lemma~11
are automatically satisfied if $\bk$ is a field. The above constructions and
results generalize well known ones for $q=-1$ $(N=2)$.\\
By starting from graded $q$-differential ideals of $\Omega_q(\cala)$ generated
in degree greater than or equal to 1, one constructs graded $q$-differential
algebras which coincide with $\cala$ in degree 0 and are generated by $\cala$.

\newpage

\section{Further results and perspectives}

Let $E$ be a $N$-differential module with $N$-differential $d$ and let us
extend the homomorphisms $[i]:H_{(m)}(E)\rightarrow H_{(m+1)}(E)$ and
$[d]:H_{(m+1)} (E)\rightarrow H_{(m)}(E),\linebreak[4] m\in \{1,\dots,N-2\}$ to
the module $H_{(\bullet)}(E)=\oplus^{N-1}_{m=1}H_{(m)}(E)$ by
setting\linebreak[4] $[i](H_{(N-1)}(E))=0$ and $[d](H_{(1)}(E))=0$. It is clear
that then, $[i]$ and $[d]$ are two $(N-1)$-differentials on $H_{(\bullet)}(E)$,
(i.e. $[i]^{N-1}=0$ and $[d]^{N-1}=0$). What is less obvious but has been
proved by Kapranov in \cite{Kapr} in that $[i]+[d]$ is also a
$(N-1)$-differential on $H_{(\bullet)}(E)$ i.e. that one has
$([i]+[d])^{N-1}=0$. In the case where $N=3$, $[i]+[d]$ is an ordinary
differential and it was also proved in \cite{Kapr} that then the differential
module $(H_{(\bullet)}(E),[i]+[d])$ is acyclic, (i.e. $\ker([i]+[d])=\im
([i]+[d]$)). It is easy to see that the latter result is equivalent to the
exactness of the hexagon $(\calh^{1,1})$ for $N=3$, i.e. it is equivalent to
Lemma 1 in the case $N=3$. In this sense, Lemma 1 appears as the generalization
to arbitrary $N\geq 3$ of this result of \cite{Kapr}.\\

It was shown in \cite{Kapr} that if $q\in \bk$ is such that Assumptions $(A_1)$
is satisfied then one can construct a tensor product for $N$-complexes in the
following manner. Let $(E_1=\oplus_nE^n_1,d_1)$ and $(E_2=\oplus_nE^n_2,d_2)$
be two $N$-complexes and let us define $d$ on $E_1\otimes E_2$ by setting
\[
d(\alpha_1\otimes \alpha_2)=d_1(\alpha_1)\otimes
\alpha_2+q^{a_1}\alpha_1\otimes d_2(\alpha_2), \ \forall \alpha_1\in
E^{a_1}_1,\ \forall \alpha_2\in E_2,
\]
one has again by induction on $n\in \mathbb N$
\[
d^n(\alpha_1\otimes \alpha_2)=\sum^n_{m=0} q^{a_1(n-m)}\left[ \begin{array}{l}
n\\ m \end{array}\right]_q d^m_1(\alpha_1)\otimes d^{n-m}_2(\alpha_2)
\]
and therefore Assumption $(A_1)$ implies $d^N(\alpha_1\otimes
\alpha_2)=d^N_1(\alpha_1)\alpha_2+\alpha_1 d^N_2(\alpha_2)=0$. Unfortunately,
as pointed out in \cite{ASi}, when $(E_1,d_1)$ and $(E_2,d_2)$ are furthermore
two graded $q$-differential algebras, $d$ fails to be a $q$-differential in
that it does not satisfy the graded $q$-Leibniz rule except for $q=1$ or
$q=-1$. This latter draw back is closely related to the fact that there is no
natural bilinear bracket on the graded $q$-derivations except in the cases
where $q=1$ or $q=-1$,\\

For $q=-1$, it was shown in \cite{Kar} that the graded commutators in
$\Omega(\cala)=\Omega_{(-1)}(\cala)$ are stable by the differential and that
the quotient complex of $\Omega(\cala)$  by the graded commutators which was
called {\sl the noncommutative de Rham complex} has a cohomology which
identifies generically with the reduced cyclic homology of $\cala$.
Unfortunately one cannot define naively along the same line a $q$-analog of the
noncommutative de Rham complex because for $q\not= 1$ and $q\not=-1$ the graded
$q$-commutators in $\Omega_q(\cala)$ are not stable by the $q$-differential.
This is again closely related to the non existence of a natural bilinear
bracket on the $q$-derivations for $q\not= 1$ and $q\not= -1$.\\
Let us finally mention a potential application of $N$-complexes in
nonassociative algebra. There is a simple-minded way to understand the usual
complexes associated to Lie algebras and to associative algebras. Indeed in
these cases the Jacobi identity for Lie algebras and the associativity
condition for associative algebras are identities which are quadratic in the
product. Therefore if one interprets the product or its dual as the first step
of a differential $d$, the above quadratic identities correspond to the first
step of $d^2=0$ (one then extends $d$ appropriately and combines it with
representations). This suggests that if one starts with a nonassociative
algebra where one has for the product identities of degree $N$ replacing
associativity, one should associate $N$-complexes to it and that a good
(co)homological theory for such an algebra would be a generalized (co)homology
of the type considered here. In particular we expect a theory with $d^3=0$ for
Jordan algebras. Work on such an approach is currently in progress.\\

It is worth noticing here that Kassel and Wambst are developing a version of
the theory of $N$-complexes which is more in the line of modern homological
algebra and which contains in particular an appropriate generalization of the
functors Tor and Ext \cite{KW}.\\

\noindent{\bf Acknowledgements.} It is a pleasure to thank warmly Mikhail
Kapranov, Max Karoubi, Christian Kassel, Richard Kerner, Peter Michor, Michel
Rausch de Traubenberg, Jim Stasheff and Marc Wambst for helpful suggestions and
comments.


\newpage

\begin{thebibliography}{999}

\bibitem[MD-V]{MD-V} M. Dubois-Violette.
 \newblock Generalized differential spaces with $d^N=0$ and the
$q$-differential calculus.
\newblock
L.P.T.H.E.-ORSAY 96/75, available from http://qcd.th.u-psud.fr; q-alg/9609012.
{\sl Czech J. Phys.\/ \bf 46} (1997) 1227-1233.

\bibitem[D-VK]{D-VK} M. Dubois-Violette, R. Kerner.
\newblock Universal $q$-differential calculus and $q$-analog of homological
algebra.
\newblock L.P.T.H.E.-ORSAY 96/48, available from http://qcd.th.u-psud.fr;
q-alg/9608026. {\sl Acta Math. Univ. Comenianae\/ \bf 65} (1996) 175-188.

\bibitem[D-VT]{D-VT} M. Dubois-Violette, I.T. Todorov.
\newblock Generalized cohomologies and the physical subspace of the $SU(2)$
WZNW model.
\newblock L.P.T.H.E.-ORSAY 97/08 and BLTP-JINR, Dubna E2-97-113,
available from http://qcd.th.u-psud.fr; hep-th/9704069. To appear in {\sl Lett.
Math. Phys.}.

\bibitem[Gre]{Gre} W.H. Greub.
\newblock Linear algebra.
\newblock Springer-Verlag 1963.

\bibitem[Kapr]{Kapr} M.M. Kapranov.
\newblock On the q-analog of homological algebra.
\newblock Preprint Cornell University 1991; q-alg/961005.

\bibitem[Kar]{Kar} M. Karoubi.
\newblock Homologie cyclique des groupes et alg\`ebres.
\newblock {\sl Acad. Sci. Paris\/ \bf 297}, S\'erie I (1983), 381-384.\\
M. Karoubi.
\newblock Homologie cyclique et $K$-th\'eorie.
\newblock {\sl Ast\'erisque\/ \bf 149} (S.M.F.), 1987.

\bibitem[KW]{KW} C. Kassel, M. Wambst.
\newblock Alg\`ebre homologique des $N$-complexes et homologie de Hochschild
aux racines de l'unit\'e.
\newblock q-alg/9705001. To appear in {\sl Pub. Res. Inst. Math. Sci. Kyoto
University.}

\bibitem[J-LL] {J-LL} J.-L. Loday.
\newblock Cyclic homology.
\newblock Springer Verlag 1992.

\bibitem[Mass] {Mass} T. Masson.
\newblock G\'eom\'etrie non commutative et applications \`a la th\'eorie des
champs.
\newblock Thesis, Orsay 1995, available from http://qcd.th.u-psud.fr.

\bibitem[May]{May} W. Mayer.
\newblock A new homology theory I,II.
\newblock {\sl Annals of Math.\/ \bf 43} (1942) 370-380 and 594-605.

\bibitem[RdT]{RdT} M. Rausch de Traubenberg.
\newblock Alg\`ebres de Clifford, supersym\'etrie et sym\'etries $\mathbb Z_n$.
Application en th\'eorie des champs.
\newblock Preprint Strasbourg LPT-02.

\bibitem[ASi]{ASi} A. Sitarz.
\newblock On the tensor product construction for $q$-differential algebras.
\newblock q-alg/9705014.

\bibitem[Weib]{Weib} C.A. Weibel.
\newblock An introduction to homological algebra.
\newblock Cambridge University Press 1994.


\end{thebibliography}
\end{document}